\documentclass{amsart}
\usepackage[all]{xy}
\usepackage{amssymb}
\usepackage{amsmath}
\usepackage{amsmath}
\swapnumbers
\newtheorem{thm}[subsection]{Theorem}
\newtheorem{prop}[subsection]{Proposition}
\newtheorem{lemma}[subsection]{Lemma}
\newtheorem{cor}[subsection]{Corollary}

\def\d{\succ}
\def\g{\prec}

\def\u{\underline }

\def\TT{T\!\!\! T}
\def\SS{\Sigma}
\def\oo{\omega}
\def\dd{\delta}
\def\DD{\Delta}

\def\aa{\alpha}

\def\t{\otimes}
\def\P{\noindent \emph{Proof.} }

\def\ue{universal enveloping }

\newcommand{\Vt}[1]{V^{\otimes #1}}

\def\HH{{\mathcal{H}}}
\def\PP{{\mathcal{P}}}
\def\RRR{{\mathcal{R}}}

\def\gg{\mathfrak{g}}
\def\cc{\gamma}

\def\Prim{\mathrm{Prim\, }}
\def\Id{\mathrm{Id }}
\def\Hom{\mathrm{Hom}}

\def\Bi{\mathop{B_{\infty}}}

\def\Im{\mathop{\rm Im}}
\def\Ker{\mathop{\rm Ker}}
\def\ue{universal enveloping }

\def\dda{2-associative }

\def\TTi{\TT_{\infty}}
\def\uu#1{{\underline #1}}
\def\DDo{\overline {\Delta}}
\def\HHo{\overline {\HH}}
\def\To{\overline {T}}

\def\arbreA{\vcenter{\xymatrix@R=3pt@C=3pt{
&& \\
&*{}\ar@{-}[ur] \ar@{-}[ul] \ar@{-}[d]     &\\
&&
}}}

\def\arbreAB{\vcenter{\xymatrix@R=2pt@C=2pt{
&&&&\\
&&&*{}\ar@{-}[ul] & \\
&&*{}\ar@{-}[uurr] \ar@{-}[uull] \ar@{-}[d]     &&\\
&&&&
}}}

\def\arbreBA{\vcenter{\xymatrix@R=2pt@C=2pt{
&&&&\\
&*{}\ar@{-}[ur] &&& \\
&&*{}\ar@{-}[uurr] \ar@{-}[uull] \ar@{-}[d]     &&\\
&&&&
}}}

\def\arbreBB{\vcenter{\xymatrix@R=2pt@C=2pt{
&&&&\\
&&&& \\
&&*{}\ar@{-}[uurr] \ar@{-}[uull] \ar@{-}[d] \ar@{-}[uu]     &&\\
&&&&
}}}

\def\arbreABC{\vcenter{\xymatrix@R=1pt@C=1pt{
&&&&&&\\
&*{}\ar@{-}[ur] &&&&& \\
&&*{}\ar@{-}[uurr] &&&&\\
&&&*{}\ar@{-}[uuurrr] \ar@{-}[uuulll] \ar@{-}[d] &&&\\
&&&&&&
}}}

\def\arbreBAC{\vcenter{\xymatrix@R=1pt@C=1pt{
&&&&&&\\
&&&*{}\ar@{-}[ul] &&& \\
&&*{}\ar@{-}[uurr] &&&&\\
&&&*{}\ar@{-}[uuurrr] \ar@{-}[uuulll] \ar@{-}[d] &&&\\
&&&&&&
}}}

\def\arbreACA{\vcenter{\xymatrix@R=1pt@C=1pt{
&&&&\\
&*{}\ar@{-}[ur] &&*{}\ar@{-}[ul] & \\
&&*{}\ar@{-}[uurr] \ar@{-}[uull] \ar@{-}[d] &&\\
&&&&
}}}

\def\arbreCAB{\vcenter{\xymatrix@R=1pt@C=1pt{
&&&&&&\\
&&&*{}\ar@{-}[ur] &&& \\
&&&&*{}\ar@{-}[uull] &&\\
&&&*{}\ar@{-}[uuurrr] \ar@{-}[uuulll] \ar@{-}[d] &&&\\
&&&&&&
}}}

\def\arbreCBA{\vcenter{\xymatrix@R=1pt@C=1pt{
&&&&&&\\
&&&&&*{}\ar@{-}[ul] & \\
&&&&*{}\ar@{-}[uull] &&\\
&&&*{}\ar@{-}[uuurrr] \ar@{-}[uuulll] \ar@{-}[d] &&&\\
&&&&&&
}}}

\def\arbreBBB{\vcenter{\xymatrix@R=1pt@C=1pt{
&&&&\\
&&&& \\
&&*{}\ar@{-}[uurr] \ar@{-}[uur] \ar@{-}[uul] \ar@{-}[uull] \ar@{-}[d]  &&\\
&&&&
}}}

\def\arbreBBA{\vcenter{\xymatrix@R=1pt@C=1pt{
&&&&\\
&&& *{}\ar@{-}[ul] & \\
&&*{}\ar@{-}[uurr] \ar@{-}[uul] \ar@{-}[uull] \ar@{-}[d]  &&\\
&&&&
}}}

\def\arbreABB{\vcenter{\xymatrix@R=1pt@C=1pt{
&&&&\\
&*{}\ar@{-}[ur] &&& \\
&&*{}\ar@{-}[uurr] \ar@{-}[uur] \ar@{-}[uull] \ar@{-}[d]  &&\\
&&&&
}}}

\def\arbreAe{\vcenter{\xymatrix@R=3pt@C=3pt{
&& \\
&*{}\ar@{-}[ur] \ar@{-}[ul]    &\\
&{*}&
}}}

\def\arbreBAe{\vcenter{\xymatrix@R=2pt@C=2pt{
&&&&\\
&&&*{}\ar@{-}[ul] & \\
&&*{}\ar@{-}[uurr] \ar@{-}[uull]     &&\\
&&{*}&&
}}}

\def\arbreABe{\vcenter{\xymatrix@R=2pt@C=2pt{
&&&&\\
&*{}\ar@{-}[ur] &&& \\
&&*{}\ar@{-}[uurr] \ar@{-}[uull]      &&\\
&&{*}&&
}}}

\def\arbreBBe{\vcenter{\xymatrix@R=2pt@C=2pt{
&&&&\\
&&&& \\
&&*{}\ar@{-}[uurr] \ar@{-}[uull]  \ar@{-}[uu]     &&\\
&&{*}&&
}}}

\def\arbreAdot{\vcenter{\xymatrix@R=3pt@C=3pt{
&& \\
&*{}\ar@{-}[ur] \ar@{-}[ul]      &\\
&\cdot&
}}}

\def\arbreABdot{\vcenter{\xymatrix@R=2pt@C=2pt{
&&&&\\
&&&*{}\ar@{-}[ul] & \\
&&*{}\ar@{-}[uurr] \ar@{-}[uull]     &&\\
&&\cdot&&
}}}

\def\arbreBAdot{\vcenter{\xymatrix@R=2pt@C=2pt{
&&&&\\
&*{}\ar@{-}[ur] &&& \\
&&*{}\ar@{-}[uurr] \ar@{-}[uull]     &&\\
&&\cdot&&
}}}

\def\arbreBBdot{\vcenter{\xymatrix@R=2pt@C=2pt{
&&&&\\
&&&& \\
&&*{}\ar@{-}[uurr] \ar@{-}[uull]  \ar@{-}[uu]     &&\\
&&\cdot&&
}}}

\def\infinitesimalun{\vcenter{\xymatrix@R=2pt@C=2pt{
&&\\
*{}\ar@{-}[u]\ar@{-}[dr]&&*{}\ar@{-}[u]\ar@{-}[dl]\\
&*{}\ar@{-}[d]&\\
&*{}\ar@{-}[dl]\ar@{-}[dr]\\
*{}\ar@{-}[d]&&*{}\ar@{-}[d]\\
&&
}}}

\def\infinitesimaldeux{\vcenter{\xymatrix@R=2pt@C=2pt{
&&&&\\
&&&*{}\ar@{-}[u]\ar@{-}[dl]\ar@{-}[dr]&\\
&&*{}\ar@{-}[d]&&*{}\ar@{-}[ddd]\\
*{}\ar@{-}[uuu]\ar@{-}[dr]&&*{}\ar@{-}[dl]&&\\
&*{}\ar@{-}[d]&&&\\
&&&&
}}}

\def\infinitesimaltrois{\vcenter{\xymatrix@R=2pt@C=2pt{
&&&&\\
&*{}\ar@{-}[dl]\ar@{-}[dr]\ar@{-}[u]&&&\\
*{}\ar@{-}[ddd]&&*{}\ar@{-}[d]&&\\
&&*{}\ar@{-}[dr]&&*{}\ar@{-}[dl]*{}\ar@{-}[uuu]\\
&&&*{}\ar@{-}[d]&\\
&&&&
}}}

\def\infinitesimalquat{\vcenter{\xymatrix@R=2pt@C=2pt{
&\ar@{-}[ddddd]&\ar@{-}[ddddd]\\
&&\\
&&\\
&&\\
&&\\
&&
}}}

\def\arbregammapq{\vcenter{\xymatrix@R=2pt@C=2pt{
&& & & & & & \\
&& *{}\ar@{-}[ul]\ar@{-}[u]\ar@{-}[ur] & & & & *{}\ar@{-}[ul]\ar@{-}[ur]& \\
\cc_{32}=&& & &*{}\ar@{-}[ull] \ar@{-}[urr] \ar@{-}[d] & & & \\
&& & & & & & \\
}}}

\begin{document}

\author[J.-L. Loday and M. Ronco]{Jean-Louis Loday \and Mar\'\i a Ronco}
\address{Institut de Recherche Math\'ematique Avanc\'ee\\
    CNRS et Universit\'e Louis Pasteur\\
    7 rue R. Descartes\\
    67084 Strasbourg Cedex, France}
\email{loday@math.u-strasbg.fr}
\urladdr{www-irma.u-strasbg.fr/{$\sim$}loday/}

\address{Depto. de Matem\`aticas, Univ. de Valpara\'\i so
                                        Avda. Gran Breta\~na 1091\\
                                        Valpara'so. Chile}
\email{maria.ronco@uv.cl}

\title{On the structure of cofree Hopf algebras}
%\alttitle{} 
\subjclass[2000]{16A24, 16W30, 17A30, 18D50, 81R60.}
\keywords{ Bialgebra, Hopf algebra, Milnor-Moore, 
Poincar\'e-Birkhoff-Witt, shuffle, 2-associative, $\Bi$-algebra, unitary infinitesimal, operad}

%\date{\today}

\begin{abstract}
We prove an analogue of the Poincar\'e-Birkhoff-Witt theorem and of the Cartier-Milnor-Moore theorem for non-cocommutative Hopf algebras. The primitive part of a cofree Hopf algebra is a $\Bi$-algebra. We construct a universal enveloping functor $U2$ from $\Bi$-algebras to $2$-associative algebras, i.e.~algebras equipped with two associative operations. We show that any cofree Hopf algebra $\HH$ is of the form $U2(\Prim \HH)$. We take advantage of the description of the free $2as$-algebra in terms of planar trees to unravel the structure of the operad $\Bi$.
\end{abstract}

\maketitle

\section*{Introduction} \label{S:int}
Any connected cocommutative Hopf algebra $\HH$ defined over a field of characteristic zero  is of the form $U(\Prim \HH )$, 
where the primitive part $\Prim \HH$ is viewed as a Lie
algebra, and
$U$ is the universal enveloping functor. This result is known as the Cartier-Milnor-Moore theorem in the literature, cf.~\cite{Ca}, \cite{MM} and \cite{Qui} appendix B. Combined with the 
Poincar\'e-Birkhoff-Witt theorem, it gives an equivalence between
the cofree cocommutative Hopf algebras and the Hopf algebras of the 
form $U(\gg)$, where $\gg$ is a Lie algebra.

Our aim is to prove a similar result without the assumption 
``cocommutative" and to get a structure theorem for cofree
Hopf algebras.
In order to achieve this goal we need to consider $\Prim \HH$ as a 
$\Bi$-algebra (cf.~1.4) instead of a Lie algebra. A
  $\Bi$-algebra is defined by $(p+q)$-ary operations for any pair of 
positive integers $(p,q)$ satisfying some relations.
This structure appears naturally in algebraic topology.
(cf.~\cite{Bau}, \cite{GJ}, \cite{Kad},  \cite{Vor}). 
The universal enveloping functor $U$ is replaced by
a functor
$U2$ from
$\Bi$-algebras to 2-associative algebras, which are vector spaces
equipped with two associative operations sharing the same unit. We prove the following
\\

\noindent {\bf Theorem.}  {\it If $\HH$ is a bialgebra over a field $K$, then 
the following are equivalent:

\noindent (a) $\HH$ is a connected \dda  bialgebra,

\noindent (b) $\HH$ is isomorphic to $U2({\Prim} {\HH})$ as a \dda bialgebra,

\noindent (c) $\HH$ is cofree among the connected coalgebras.}
%\vspace*{2pt}

Hence, as a consequence, we get a structure theorem for cofree Hopf 
algebras: any cofree Hopf algebra is of the form $U2(R)$ where $R$ is a $\Bi$-algebra. 
The notion of \dda bialgebra occuring in the theorem is as follows. A 
\dda bialgebra $\HH$ is a vector space equipped with two associative operations
denoted
$*$ and $\cdot$ and a cooperation $\DD$. We suppose that $(\HH, *, 
\DD)$ is a bialgebra (in the classical sense), and that
$(\HH, \cdot, \DD)$ is a \emph {unital infinitesimal} bialgebra. The difference 
with the classical notion is in the \emph {compatibility relation} between
the product and the coproduct which, in the unital infinitesimal case, 
is
$$ \DD(x\cdot y) = (x\t 1) \cdot \DD(y) + \DD(x) \cdot (1\t y) - x\t y 
\ .$$
The following structure theorem for unital infinitesimal bialgebras is the  key result in the proof of the main theorem.
\\

\noindent {\bf Theorem.} \emph{
Any connected unital infinitesimal bialgebra is isomorphic to the tensor algebra (i.e.~non-commutative polynomials) equipped with the deconcatenation coproduct.}

Observe that, in our main theorem, (a) $\Rightarrow$ (b) is the analogue of 
the Cartier-Milnor-Moore theorem (cf.~\cite{Ca}, \cite {MM}), that  (b) $\Rightarrow$
(c) is the analogue of the Poincar\'e-Birkhoff-Witt theorem, and that  
(a) $\Rightarrow$ (c) is the analogue of a theorem of Leray in the
cocommutative case.

The universal enveloping $2as$-algebra $U2(R)$, where $R$ is a $\Bi$-algebra, is a quotient of the free $2as$-algebra $2as(V)$ for $V=R$. So it is important to know an explicit description of $2as(V)$ for any vector space $V$. We show that, as an associative algebra for the product $*$, $2as(V)$ is the non-commutative polynomial algebra (i.e.~tensor algebra) over the planar rooted trees. We describe the other product $\cdot$ and the coproduct $\DD$ in terms of trees. This description is very similar to the description of the Hopf algebra of (non-planar) rooted trees given by Connes and Kreimer (cf. \cite{CK}). Here the rooted trees are replaced by the planar rooted trees, and the polynomials by the non-commutative polynomials.

  The free Lie algebra is a complicated object which is mainly studied 
through its embedding in the free associative algebra, that is  its
identification with the primitive part of the tensor algebra. The free 
$\Bi$-algebra is, a priori,  an even more complicated object (it is
generated by $k$-ary operations for all $k$). Our main theorem permits us to identify it with
the primitive part of the free  \dda algebra. Since we can describe explicitly the free  \dda algebra in terms of trees, we can prove the following
\\

\noindent {\bf Corollary.} \emph{The operad $\Bi$ of $\Bi$-algebras is such that $\Bi(n)= K[\TT_n]\t K[S_n]$, where $\TT_n$ is the set of planar rooted trees with $n$ leaves, and $S_n$ is the symmetric group.}
\\

The content of this paper is as follows. In section 2 we recall the classical 
notion of bialgebras and $\Bi$-algebras.
In section 3, which can be read independently of section 2,  we recall the notion of unital infinitesimal bialgebra and we prove the key result which says that there is only
one kind of connected unital infinitesimal bialgebra. In
section 4 we introduce the notion of \dda algebra and
\dda bialgebra, and we construct the universal enveloping functor $U2$. In section 5 we state and 
prove the main theorem. Our proof is based on the study of
the free objects and minimizes combinatorics computation. In section 6 we give an explicit description of the free 2-associative algebra in terms of planar trees.   Then we 
prove that the  operad of 2-associative algebras is a Koszul operad in
the sense of Ginzburg and Kapranov, and we describe the chain complex 
which computes the homology of a
\dda algebra (a Hochschild type complex).  In section 7 we unravel the structure of bialgebra of the free $2as$-algebra. We show that it has a description analogous to the Connes-Kreimer Hopf algebra. We prove that it is self-dual. We deduce from our main theorem and the previous section an explicit description of  the free $\Bi$-algebra. In section 8 we state without proofs
a variation of our main result in terms of ``dipterous" 
algebras. It has the advantage of taking care of the cofree Hopf algebras whose primitive part is in fact a brace algebra. Then, in section 9,  we compare several results of the same kind.

The main result of this paper was announced without proof in \cite{LR2}.
\\

\noindent {\bf Notation.} In this paper $K$ is a field and all vector spaces 
are over $K$. Its unit is denoted $1_K$ or just 1. The vector space spanned by the elements of a set $X$ is denoted $K[X]$. The tensor product 
of vector spaces over $K$
is denoted by $\t$. The tensor product of $n$ copies of the space $V$ is
denoted $\Vt n$. For $v_i\in V$ the element $v_1\t \cdots \t v_n$ of 
$\Vt n$ is denoted by  $(v_1, \ldots , v_n)$ or
simply by $v_1 \ldots v_n$.
A linear map $\Vt n \to V$ is called an {\it $n$-ary operation} on $V$ 
and a linear map $ V\to \Vt n$ is called an {\it
$n$-ary cooperation} on $V$.

 \section{Hopf algebra and $\Bi$-algebra}\label{S:Hopf}
We recall the definition of Hopf algebra,
  tensor algebra, tensor coalgebra and the definition of
$\Bi$-algebra together with their relationship.

\subsection{Hopf algebra}\label{Hopf} By definition a {\it bialgebra} $(\HH, *, \DD, u, c)$ is a vector space 
$\HH$ equipped with an associative product
$*:\HH\t \HH\to \HH$, a unit  $u:K\to \HH$, and a coassociative 
coproduct $\DD:\HH\to \HH\t \HH$, a counit
$c:\HH\to K$  such that $*$ and $u$ are morphisms of coalgebras or, 
equivalently, $\DD$ and $c$ are
morphisms of algebras.

We will use the notation $\HHo:= \Ker c$, so that $\HH = K 1 
\oplus \HHo$, and the notation $\DD(x) = x \t 1  +1\t x
+  \DDo (x)$ for
$x\in \HHo$. In general we will omit $u$ and $c$ in the notation, the 
unit of $\HH$, that is $u(1_K)$, being denoted
by $1$.

Given two linear maps $f,g : \HH \to \HH$ their {\it convolution} is, 
by definition, the composite
  $f\star g := *\circ (f\t g)\circ \DD$. The convolution is associative 
with unit the map $u\circ c$. By definition an {\it antipode}
  for $\HH$ is a linear map $S:\HH \to \HH$ such that $S\star \Id = 
u\circ c = \Id \star S$.

By definition a {\it Hopf algebra} is a   bialgebra equipped with an 
antipode. It is well-known that for a
conilpotent bialgebra such an antipode automatically exists (see below).

An element $x$ of $\HH$ is said to be {\it primitive} if $\DD(x) = x\t 
1 + 1 \t x$ or, equivalently, $\DDo  (x) = 0$. The
space of primitive elements of $\HH$ is denoted $\Prim \HH$. It is 
known to be a sub-Lie algebra of $\HH$ for the Lie
structure given by $[x,y] := x*y - y*x$.

\subsection{$n$-ary (co)-operations, connectedness}\label{connectedness}
For an associative algebra there is (essentially) only one $(n+1)$-ary 
operation. It is given by $*^{n}(x_0\cdots x_{n}) =
x_0*\cdots * x_{n}$.

Dually $\DDo $ determines an $(n+1)$-ary cooperation $\DDo ^n:\HH \to 
\HH^{\t n+1}$ given by $\DDo ^0 = \Id,\
\DDo ^1 = \DDo  $, and
$\DDo ^n = (\DDo \t \Id\t  \cdots \t \Id )\circ \DDo ^{n-1}$.

Following Quillen (cf.~\cite{Qui}, p. 282) we say that a bialgebra $\HH$ is 
{\it connected} if $\HH = \bigcup _{r\geq 0} F_r\HH$
where $F_r\HH$ is the coradical filtration of $\HH$ defined recursively by the 
formulas
\begin{eqnarray*}
F_0\HH &:=& K 1, \\
F_r\HH & :=& \{ x\in \HH \mid  \DDo  (x) \in F_{r-1}\HH \t F_{r-1}\HH \}\ .
\end{eqnarray*}

Observe that we use only $\DD$ and 1 to define connectedness.

If $\HH$ is connected, then $\HH$ is conilpotent, that is for any 
element
$x\in \HHo $ there exists $n$ such that $\DDo  ^n(x)=0$. In this case 
the antipode $S$ is given by
\begin{displaymath}
S(x) := \sum_{n\geq 0} (-1)^{n+1} *^n \circ {\DDo }^n (x)\ .
\end{displaymath}

Therefore a connected bialgebra is equivalent to a connected Hopf 
algebra. In this paper we will mainly work with
connected bialgebras.

\subsection{Tensor algebra and tensor coalgebra}
By definition the {\it tensor algebra} $T(V)$ over the vector space $V$ 
is the {\it tensor module}
$$T(V) = K \oplus V \oplus \Vt 2 \oplus \cdots \oplus\Vt n \oplus 
\cdots $$
equipped with the associative product called {\it concatenation}:
$$v_1\cdots v_i \t v_{i+1} \cdots v_n \mapsto v_1 \cdots v_i v_{i+1}\cdots 
v_n \ . $$
It is known to be the free associative algebra over $V$. It is an 
augmented algebra, whose augmentation ideal is
denoted by
$\To(V)$.  It is well-known that $T(V)$ is a connected Hopf algebra for 
the coproduct defined by the shuffles. Its
primitive part is the free Lie algebra on $V$.

Dually the  {\it tensor coalgebra} $T^c(V)$ over the vector space $V$ 
is the tensor module  (as above)
equipped with the coassociative coproduct $\DD$ called {\it 
deconcatenation}:
\begin{displaymath}
\DD( v_1 \cdots v_n ) = \sum_{i=0}^{i=n}v_1\cdots v_i \t  v_{i+1} 
\cdots v_n \ . 
\end{displaymath}
Observe that it is a connected coalgebra and that its primitive part is 
$V$. It satisfies the following universal condition. Given a connected coalgebra $C$ and a linear map 
 $\phi : C \to V$ such that $\phi(1)=0$, there exists a unique  coalgebra map $\overline { \phi} : C \to T^c(V)$ which extends $\phi$.
 Explicitly the $k$th component of   $\overline \phi (c)$ is given by 
$\phi^{\t k}\circ {\DDo}^{k-1} (c)$ for $k\geq 2$. Therefore $T^c(V)$ is \emph {cofree} (among connected coalgebras). 

We will say that a bialgebra $\HH$ is {\it cofree} if, as a 
coalgebra, it is isomorphic to $T^c(\Prim \HH)$. Such an isomorphism is completely determined by a projection 
$\HH  \to \Prim \HH$ inducing the identity on $\Prim \HH$ and sending 1 to 0. Indeed, since
$\HH$ is cofree, this projection extends uniquely as coalgebra morphism 
$\HH \to T^c(\Prim \HH)$. A bialgebra which is
cofree is, by definition, connected and therefore is a Hopf algebra.

The tensor coalgebra is known to be a Hopf algebra for the product 
induced by the shuffles:
$$v_1\cdots v_p\, {\scriptstyle \sqcup\!\sqcup} \, v_{p+1} \cdots v_{p+q} := \sum 
v_{i_1} \cdots v_{i_{p+q}}$$
where the sum is extended to all permutations $({i_1}, \cdots 
,i_{p+q})$ of $(1, \cdots , p+q)$ which are
$(p,q)$-shuffles, i.e.~such that $(1, \cdots , p)$ appear in this order 
and  $(p+1, \cdots , p+q)$  appear in this
order. We will call it the {\it shuffle algebra} and denote it by 
$T^{sh}(V)$.

\subsection{$\Bi$-algebra}\label{Binfini}
Let $\HH$ be a cofree bialgebra and put $V:= \Prim \HH$, so we have 
$\HH \cong T^c(V)$ as a coalgebra.
Transporting the algebra structure of $\HH$ under this isomorphism 
gives a coalgebra homomorphism
$$* : T^c(V)\t T^c(V) \to T^c(V)\ .$$
Since $T^c(V)$ is cofree the map $*$ is  completely determined by its 
value in $V$ (degree one component of
$T^c(V)$), that is by maps
$$M_{pq}: \Vt p \t \Vt q \to V\ ,\  p\geq 0, q\geq 0.$$
From the unitality and counitality property of $\HH$ we deduce that
$$M_{00}= 0, \ M_{10}= \Id_V = M_{01}, \ {\rm and }\ M_{n0}=0=M_{0n}\ 
{\rm for} \ n\geq 2.$$
For any  set of indices $(\u i,  \u j) := (i_1, \ldots, i_k ; j_1,
\ldots , j_k)$ such that $i_1+ \cdots+ i_k = p$ and  $j_1+ \cdots+j_k = 
q$ we denote by
$$M_{i_1j_1}M_{i_2j_2}\ldots M_{i_kj_k}:  \Vt p \t \Vt q \to \Vt k$$
  the map which sends
$(u_1u_2\ldots u_p,v_1v_2\ldots v_q)\in  \Vt p \t \Vt q$ to
\begin{displaymath} \begin{array}{l}
M_{i_1j_1}(u_1\ldots u_{i_1},v_1\ldots v_{j_1})M_{i_2j_2}(u_{i_1+1}\ldots u_{i_1+i_2},v_{j_1+1}\ldots 
v_{j_1+j_2})\ldots \qquad\qquad \\
\hfill \ldots M_{i_kj_k}( \ldots u_p, \ldots v_q)\in  \Vt k\ .
\end{array} \end{displaymath}
The associative operation on $T^c(V)$ is recovered from the $(p+q)$-ary 
operations $M_{pq}$ by the formula
\begin{equation}\label{Mpq} u_1\cdots u_p * v_1 \cdots v_q = \sum_{k\geq 1} \big(\sum_{(\u i, \u 
j)} M_{i_1j_1}\ldots M_{i_kj_k}
(u_1\ldots u_p,v_1\ldots v_q)\big)\ ,
\end{equation}
where the sum is extended over all sets of indices  $(\u i, \u 
j) := (i_1, \ldots, i_k ; j_1,
\ldots , j_k)$ such that $i_1+ \cdots+ i_k = p$ and  $j_1+ \cdots+j_k = 
q$.
Of course the $k$th component is in $V^{\t k}$. We adopt the notation
\begin{equation}\label{Mkij} M^k_{(\u i , \u j )}:=M_{i_1j_1}\ldots 
M_{i_kj_k}\ . 
\end{equation}
   Observe that the component
in $V$, i.e. for $k=1$, is $M_{pq}(u_1 \cdots u_p,
v_1\cdots v_q)$. The last non-trivial component is in $V^{\t p+q}$ and 
is
made of all the $(p,q)$-shuffles of $u_1 \cdots
u_p v_1 \cdots v_q$. For instance in low dimension one gets
\begin{eqnarray*}
u*v &=& M_{11}(u, v) + uv + vu \ ,\\
uv*w &=&  M_{21}(uv , w) + u M_{11}(v, w) +  M_{11}(u, w)v + uvw + uwv+
wuv\ , \\
&=& M_{21}(uv,w) + (u*w)v + u(v*w) - uwv\ , \\
u*vw &=&  M_{12}(u, vw) + M_{11}(u, v)w +  vM_{11}(u, w) + uvw + vuw+
vwu\ , \\
  &=& M_{12}(u, vw) + (u*v)w +  v(u*w) - vuw\ .
\end{eqnarray*}
Associativity of $*$ implies some relations among the operations 
$M_{pq}$. In fact we can write one such relation
$\RRR_{ijk}$ for any triple $(i,j,k)$ of positive integers by writing
$$(u_1\cdots u_i * v_1 \cdots v_j)*w_1\cdots w_k = u_1\cdots u_i *( v_1 
\cdots v_j*w_1\cdots w_k)$$
in terms of $M_{pq}$ and equating the components in $V$. It comes:
$$\sum_{1\leq l \leq i+j} M_{ l  k} \circ (M^l_{(\u i ,\u j)}\t \Id^{\t k} ) 
=
\sum_{1\leq m \leq j+k} M_{i m} \circ (\Id^{\t i} \t M^m_{(\u j, \u k)})\ 
.\eqno (\RRR_{ijk})$$
The first nontrivial relation, obtained by writing 
$(u*v)*w = u*(v*w)$,  reads:
$$M_{21}(uv+vu, w)+ M_{11}(M_{11}(u,v),w ) =  M_{12}(u, 
vw+wv)+M_{11}(u,M_{11}(v,w ))
\ .\eqno (\RRR_{111})$$

\subsection{Definition}\label{Def:Binfini}{\rm A {\it $\Bi$-algebra} (cf.~\cite{Bau}, \cite{GJ}, \cite{Vor}) is a 
vector space $R$ equipped with operations
$$M_{pq}: R^{\t p} \t R^{\t q} \to R\ ,\  p\geq 0, q\geq 0$$
satisfying
$$M_{00}= 0, \ M_{10}= \Id_R = M_{01}, \ {\rm and }\ M_{n0}=0=M_{0n}\ 
{\rm for} \ n\geq 2,$$
and  the relations $\RRR_{ijk}$ for any triple $(i,j,k)$ of positive 
integers.}

There are obvious notion of morphism, ideal and free object for 
$\Bi$-algebras. The free $\Bi$-algebra over the
vector space $V$ is denoted $\Bi(V)$.

In the literature a $\Bi$-algebra is often equipped by definition with a grading and 
a differential satisfying  some relations
with the other operations (cf.~\cite {GJ},\cite{Kad}, \cite{Vor}). Here we mainly consider non-graded objects 
with $0$ differential, whence our terminology. When a grading and a 
differential will be needed we will say {\it differential graded 
$\Bi$-algebra}, cf. ~\ref{differential}.

From the previous discussion it is clear that we have the following result.

\begin{prop} \label{Tfc} Any $\Bi$-algebra $R$ defines a cofree Hopf algebra\\ $(T^c(R), *, \DD)$, 
where $\DD$ is the deconcatenation and $*$ is given by formula (\ref{Mpq}).\hfill $\square$
\end{prop}

\subsection{Examples}\label{examples}
{\rm (a) If $M_{pq}=0$ for all $(p,q)$ different from $(0,1)$ and
$(1,0)$, then $R$ is just a vector space $V$
and  $(T^c(R), *, \DD)$ is the shuffle algebra $T^{sh}(V)$.
\\

\noindent (b) If $M_{pq}=0$ for all $(p,q)$ different from
$(0,1)$, $(1,0)$ and $(1,1)$, then $M_{11}$ is an associative operation 
by $\RRR_{111}$. Hence $R$ is simply  an associative algebra (possibly without unit).
The Hopf algebra $(T^c(R), *, \DD)$ is called the {\it quasi-shuffle 
algebra} over $R$. The product $*$ is completely determined by the 
inductive relation
$$a\oo * b \theta = M_{11}(a,b)(\oo * \theta) + a(\oo * b \theta) + b( 
a \oo * \theta)\ ,$$
where $a,b \in R$ and $\oo , \theta$ are tensors.
\\

\noindent (c) If $M_{pq}=0$ for all $(p,q)$ such that $p\geq 2$, then we get a 
{\it brace algebra}, cf.~ for instance \cite{Ger}, \cite{Kad}, \cite{R2}, \cite{Vor}.
\\

\noindent (d) A \emph{prop} is a family of $S_n\times S_m^{op}$-modules $\PP(n,m)$ equipped with a composition
\begin{displaymath}
\PP(n_1, m_1)\t \cdots \t \PP(n_p, m_p)\t \PP(r_1, s_1)\t \cdots \t \PP(r_q, s_q)
\stackrel{\cc}{\to}
\PP(n_1+\cdots n_p, s_1+\cdots s_q)
\end{displaymath}
for $m_1+\cdots +m_p=r_1+\cdots +r_q$, which is compatible with the action of the symmetric groups and is associative in an obvious sense. 
It gives rise to a $\Bi$-algebra structure on $R= \bigoplus_{n,m} \PP(n,m)$ as follows. The map $M_{p,q}$ is $\cc$ on $\PP(n_1, m_1)\t \cdots 
\t \PP(n_p, m_p)\t \PP(r_1, s_1)\t \cdots \t \PP(r_q, s_q)$ if $m_1+\cdots m_p=r_1+\cdots r_q$ and $0$ otherwise.}

\subsection{Remark} Since $R$ is the primitive part of the Hopf 
algebra $T^c(R)$, it is a Lie algebra. Using formula
  $\RRR_{111}$ one sees that the Lie bracket on $R$ is given by $[r,s] = 
r*s - s*r = M_{11}(r,s)- M_{11}(s,r)$. Hence $M_{11}$ is
a Lie-admissible operation. If $R$ is a brace algebra, then $M_{11}$ is a pre-Lie operation
since the associator of $M_{11}$ is symmetric in the last two variables.

\subsection{Relationship with deformation theory}{} A $\Bi$-structure on 
the vector space $V$ is equivalent to a
deformation of the shuffle algebra $T^{sh}(V)$: make the product $*$ 
into a product in $T^{sh}(V)[[h]]$ by taking the
element in
$V^{\t p+q-i}$ as coefficient of $h^i$.

\section{Unital infinitesimal bialgebra}\label{S:infinitesimal}
We recall from \cite{Lod2}  the notion of unital 
infinitesimal bialgebra and we prove a structure theorem.

\subsection{Definition}\label{Def:infinitesimal}
 {\rm A {\it unital infinitesimal bialgebra} $(\HH, 
\cdot, \DD)$ is a vector space $\HH$ equipped with a unital associative 
product $\cdot$ and a counital coassociative coproduct $\DD$ which are 
related by the  {\it unital infinitesimal relation}:
\begin{equation}\label{ui}
 \DD(x\cdot y) = (x\t 1) \cdot \DD(y) + \DD(x) \cdot (1\t y) - x\t y 
\ .
\end{equation}
Here the product $\cdot$ on $\HH \t \HH$ is given by
\[(x\t y)\cdot (x'\t y') := x\cdot x' \t y\cdot  y'\ .\]
Pictorially this relation reads:

$$\infinitesimalun = \infinitesimaldeux + \infinitesimaltrois - \infinitesimalquat$$

Equivalently the relation verified by the reduced coproduct $\DDo$ is:
\begin{equation}\label{reduced infinitesimal} \DDo(x\cdot y) = (x\t 1) \cdot \DDo(y) + \DDo(x) \cdot (1\t y) + x\t y 
\ .
\end{equation}}

\subsection{Example} Let $K[x]$ be the polynomial algebra in $x$. The map 
$\DD$ defined by $\DD(x^n)= \sum _{p=0}^{n} x^p
\t x^{n-p}$ satisfies the unital infinitesimal relation.

The unital infinitesimal relation differs from the 
infinitesimal relation used by S.~Joni and G.-C.~Rota in \cite{JR} (see also \cite{Agu}) by the presence
  of the term $-x\t y$. From our  relation it comes $\DD(1) = 1 \t 1$. Recall that the 
notion of connectedness given in \ref{connectedness} uses only $\DD$ and 1, and
so is applicable here.

\subsection{Proposition-Notation} \emph{The tensor module over $V$ 
equipped with the concatenation product $\cdot$ and the
deconcatenation coproduct
  $\DD$ is a unital infinitesimal bialgebra denoted $T^{fc}(V)$.}

\P  Let us compute $\DD (x\cdot y)$ for $x=u_1\ldots u_p $ and 
$y= u_{p+1} \ldots u_n$ :
\begin{eqnarray*}
\DD(x\cdot y)& =& \DD(u_1\ldots u_n )\\
&=& \sum_{i=0}^{n}u_1\ldots u_i \t u_{i+1}\ldots u_n  \\
&=&  \sum_{i=0}^{p}u_1\ldots u_i \t u_{i+1}\ldots u_n - u_1\ldots u_p \t 
u_{p+1}\ldots u_n\\
&&\qquad  + \sum_{i=p}^{n}u_1\ldots u_i \t u_{i+1}\ldots u_n \\
&=&  \DD(u_1\ldots u_p)\cdot (1 \t u_{p+1}\ldots u_n) - u_1\ldots u_p \t 
u_{p+1}\ldots u_n \\
&&\qquad + (u_1\ldots u_p \t 1)\cdot \DD(u_{p+1}\ldots u_n) \\
&=& \DD (x) \cdot (1\t y) - x\t y + (x\t 1)\cdot \DD (y) \ .
\end{eqnarray*}
\hfill $\square$

\subsection{Convolution product} For any bialgebra  $(\HH, \nu, \DD)$ (either classical or unitary infinitesimal) the {\it convolution product} $\star$ on $\Hom_K(\HH, \HH)$ is defined as follows. For $f$ and $g\in \Hom_K(\HH, \HH)$ one puts
\[ f\star g := \nu \circ (f\t g)\circ \DD\ .\]
From the associativity of $\nu$ and the coassociativity of $\DD$ it follows that $\star $ is associative. It is also easy to check that $uc$ is a unit for $\star$.

\begin{prop}\label{idempotent} Let $(\HH, \nu, \DD)$ be  a connected 
unital infinitesimal bialgebra. The linear  operator
  $e:\HH \to \HH$ defined by
\[ e:= J - J\star J + \cdots + (-1)^{n-1}J^{\star n}+\cdots \]
where $J:= \Id -uc$,  has the following properties:
\\

a) $\Im e = \Prim \HH$,

b) for any $x,y\in \HHo$ one has $e(\nu(x,y))=0$,

c) $e$ is an idempotent,

d) for $\HH = (T(V), \cdot, \DD)$, where $\cdot$ is the concatenation 
and $\DD$ the deconcatenation, $e$ is the
identity on $V$ and $0$ on the other components.
\end{prop}

\P We adopt the following notation for this proof:
$\nu(x,y) := x\cdot y$, $\Id = \Id_{\HHo}$, 
and, for any $x\in \HHo$, 
\begin{equation}\label{Deltabar}
\DDo(x) := x_{(1)}\t x_{(2)}\ .
\end{equation}

So we omit the summation symbol in the formula for the reduced comultiplication.

First, we observe that, on ${\HHo}$,  $e$ can be written
\begin{equation}\label{eseries} e:=\sum _{r\geq 0}(-1)^{r} \nu^{ r}\circ \DDo^r = \Id - \nu \circ 
\DDo + \nu^{2}\circ \DDo^2 +\cdots
\end{equation}
and so satisfies the following equality:
\begin{equation}\label{eDelta}e=\Id -\nu\circ (\Id\t e)\circ {\DDo}\ .
\end{equation}
which can be written
\begin{equation}\label{eformula}e(x)= x - x_{(1)} \cdot e(x_{(2)})\ .
\end{equation}

\noindent a) We proceed by induction on the filtration-degree of $x\in \HHo$. Recall that if $x\in F_n\HHo$, then 
$x_{(1)}$ and $x_{(2)}$ are in $ F_{n-1}\HHo$.

If $x\in F_1\HHo= \Prim \HH$, then $\DDo(x) =0$ and therefore $e(x)=x$ by (\ref{eformula}). Let us now suppose that
$e(y)\in
\Prim
\HH$ for any $y\in  F_{n-1}\HHo$, and let $x\in  F_{n}\HHo$. We compute
\begin{eqnarray*}
\DD(e(x))&=&\DD (x)-\DD \circ \nu \circ (\Id\t e)\circ {\DDo }(x) \qquad \hbox{ by (\ref{eDelta})} \\
&=&x\t 1+1\t x+x_{(1)}\t x_{(2)}-\DD (x_{(1)}\cdot e(x_{(2)}))\qquad \hbox{ by (\ref{Deltabar})} \\
&=&x\t 1 +1\t x +x_{(1)}\t x_{(2)}-\DD (x_{(1)})\cdot (1\t e(x_{(2)}))-x_{(1)}\cdot e(x_{(2)})\t 1\\
 &&\qquad  \hbox{ by (\ref{ui}) and induction}\\
&=&e(x)\t 1 +1\t x +x_{(1)}\t x_{(2)}-x_{(1)}\t e(x_{(2)})-1\t x_{(1)}\cdot
e(x_{(2)})\\
&&-(\Id\t \nu )\circ (\Id\t \Id\t e)\circ {\DDo }^2(x)\qquad  \hbox{ by (\ref{eformula})} \\
&=&e(x)\t 1 +1\t e(x) +(\Id\t (\Id -e))\circ {\DDo }(x)-(\Id\t (\Id -e))\circ {\DDo }(x)\\
&&\qquad  \hbox{ by (\ref{eDelta})} \\ 
&=&e(x)\t 1 +1\t e(x)\ ,
\end{eqnarray*}
which proves that $e(x)$ is primitive. If $x\in \Prim \HH$, then $\DDo (x) = 0$ and therefore $e(x)=x$ by 
(\ref{eDelta}).
\\

\noindent b) We proceed by induction on the sum of the filtration-degrees of $x$ and $y$. If $x$ and $y$ are both
primitive, then 
\[ \DDo(x\cdot y) = x\t y \quad \hbox{ and }\quad  \DDo^r(x\cdot y) = 0 \hbox{ for } r\geq 2. \]
Therefore we get
\[ e(x\cdot y)=x\cdot y -\nu\circ \DDo(x\cdot y) =x\cdot y -x\cdot y =0. \]
We now suppose that the formula holds when the sum of the filtration-degrees is strictly less than 
the sum of the filtration-degrees of $x$ and $y$. We have 
\begin{eqnarray*}
e(x\cdot y) &=&x\cdot y - \nu \circ (\Id \t e)\circ {\DDo }(x\cdot y)\hfill  \hbox{ by (\ref{eformula})} \\
&=&x\cdot y -x\cdot y_{(1)}\cdot e(y_{(2)})-x_{(1)}\cdot e(x_{(2)}\cdot y) -x\cdot e(y)\hfill  \hbox{ by (\ref{ui})
and (\ref{Deltabar})}\\
&=&x\cdot y -x\cdot y_{(1)}\cdot e(y_{(2)})-x\cdot e(y)\hfill   \hbox{ by induction} \\
&=&x\cdot e(y)-x\cdot e(y)=0 \hfill  \hbox{ by (\ref{eDelta})} 
\end{eqnarray*}
and the proof is completed.
\\

\noindent (c) Follows immediately from (a), since $\DDo(x)= 0$ when $x$ is 
primitive.
\\

\noindent (d) This assertion is immediate by direct inspection. \hfill $\square$

\begin{thm}\label{thm:infinitesimal}
Any connected unital infinitesimal  bialgebra $\HH$ is isomorphic to 
$T^{fc}(\Prim\HH):=(T(\Prim\HH), \nu, \DD)$, where $\nu=$ concatenation and $\DD=$ deconcatenation.
\end{thm}

\P Let $V:= \Prim \HH$. Define a bialgebra morphism $G: \HHo \to \To(V)$ by the formula 
\[ G(x):= \sum_{n\geq 1} e^{\t n} \circ \DDo^{n-1}(x). \]
Define $F:\To(V) \to \HHo$ by  $F(v_1\ldots v_n):= v_1\cdot \ldots \cdot v_n$ for $n\geq 1$. The composite 
$F\circ G$ is equal to  $\sum_{n\geq 1} e^{\star n} $ since 
\[  F\circ G =  \sum_{n\geq 1} \nu^{n-1}\circ e^{\t n} \circ \DDo^{n-1} =\sum_{n\geq 1} e^{\star n} .\]
The two series  $g(t):= t-t^2+t^3-\ldots ={t\over {1+t}} $ and $f(t):= t+t^2+t^3+\ldots={t\over {1-t}} $ are inverse to each other for composition: $(f\circ g)(t)=t$. We can apply these series to elements of $\Hom_K(\HH, \HH)$ which send 1 to 0 by using $\star$ for multiplication. We get $e= g^{\star}(J)$ (cf.~Proposition \ref{idempotent}) and
\[ F\circ G = \sum_{n\geq 1} e^{\star n} =  f^{\star}(e) =  f^{\star}\circ g^{\star}(J) = 
 (f\circ g)^{\star}(J) = \Id^{\star}(J) = J. \]
On the other hand one has 
\begin{eqnarray*}
G\circ F(v_1\ldots v_n) &=& G(v_1\cdot \ldots \cdot v_n) \\
& = & e^{\t n} \circ \DDo^{n-1}(v_1\cdot \ldots \cdot v_n)\\
&=&  e(v_1)\t \ldots \t e(v_n) \\
& = &  v_1\t \ldots \t v_n\ .
\end{eqnarray*}
We have used \ref{idempotent}.b  in this computation. We have shown that $\HHo$ and $\To(V)$ are isomorphic. \hfill $\square$

\subsection{Remark} {\rm The idempotent $e$ is the geometric series (for 
convolution) applied to the map $J= \Id - uc$. It is the analogue
of the first Eulerian idempotent, which, in the classical case, is defined as 
the logarithm series applied to $J$. The geometric
series was already used in a similar context in \cite{R1}. Observe that in our case the characteristic zero hypothesis is 
not needed since the geometric series has no denominators.}

Theorem \ref{thm:infinitesimal} is similar to the Hopf-Borel theorem, cf. \cite{Bor},  which states, in the non-graded case,  that any 
 connected commutative cocommutative Hopf algebra $\HH$ is
isomorphic to the symmetric algebra $S(\Prim \HH)$ (in characteristic zero).

\section{2-associative algebra and \dda bialgebra}\label{S:3}
In this section we introduce the algebras with two associative 
operations, that we call 2-associative algebras.
Then we study the \dda algebras equipped with a cooperation.

\subsection{Definition}\label{def:2as} {\rm A {\it  2-associative algebra} over $K$ is a 
vector space $A$ equipped with two associative operations
$(x,y)\mapsto x*y$ and $(x,y)\mapsto x\cdot y$. A \dda algebra is said 
to be {\it unital} if there is an element
1 which is a unit for both operations. Unless otherwise stated we suppose 
that the \dda algebras are unital.}

Observe that the definition of  a 2-associative
object makes sense in any monoidal category.  So we can define a notion 
of {\it 2-associative  monoid}, of {\it \dda  group}, of {\it \dda monoidal category}, of {\it \dda operad}, 
etc.

The {\it free \dda  algebra}  over the vector space $V$ is the \dda 
algebra $2as(V)$ such that  any map from $V$ to a \dda algebra 
$A$ has
a natural extension as a \dda morphism $2as(V)\to A$. In other words 
the functor $2as(-)$ is left adjoint to the
forgetful functor from \dda algebras  to vector spaces. It is clear 
that $2as(V)$ is graded and of the form $2as(V)
= \oplus_{n\geq 0} 2as_n \t V^{\t n}$. More information on the explicit 
structure of $2as(V)$, that is on $2as_n$,  is given in section 6.

\subsection{Tensor product of \dda  algebras} Given two \dda  algebras 
$A$ and $B$  we define their {\it tensor product} as the
\dda algebra  $A\t B$ equipped with the two products
$$(a\t b) *  (a'\t b') := a*a' \t b* b',$$
$$(a\t b)\cdot (a'\t b') := a\cdot a' \t b\cdot  b'\ .$$
The unit of $A\t B$ is $1\t 1$.

If $f_i : A_i \to A_i' $ (for $i=1,2$) is a morphism of \dda  algebras, 
then obviously $f_1\t f_2 : A_1\t A_2 \to
  A_1'\t A_2'$ is a morphism of \dda algebras.

\subsection{On \dda algebras and $\Bi$-algebras}\label{2asAndB} Let $(A,*, \cdot)$ be 
a 2-associative algebra.
  We define
$(p+q)$-ary operations
$M_{pq}: A^{\t p+q} \to A$, for $p\geq 0, q\geq 0$ by induction as 
follows:
\[ M_{00}= 0, \ M_{10}= \Id_A = M_{01}, \ {\rm and }\ M_{n0}=0=M_{0n}\ 
{\rm for} \ n\geq 2,\]
and
\begin{eqnarray}\label{Mpqinductif}
M_{pq}(u_1 \ldots u_p, v_1\ldots v_q)
 &:=& \big(u_1 \cdot u_2\cdot  \ldots  \cdot u_p\big) * \big(v_1 \cdot 
v_2\cdot\ldots \cdot v_q\big)\nonumber\\
&-&\sum_{k\geq 2}\sum_{(\u i ,\u j)} \big(M_{i_1j_1}
\cdot\ldots\cdot
M_{i_k j_k}(u_1\ldots u_p, v_1 \ldots v_q)\big)
\end{eqnarray}
where the second sum (for which
$k\geq 2$ is fixed) is extended to all the sets of indices $(\u i,\u j)$
  such that $i_1+ \cdots+ i_k = p$ and  $j_1+ \cdots+j_k = q$.

For instance (cf.~\ref{Binfini}):
\begin{eqnarray*}
M_{11}(u, v)&=& u*v -  u\cdot v - v\cdot u\ ,\\
M_{21}(uv , w)&=& (u\cdot v)*w   - u\cdot M_{11}(v, w) -   M_{11}(u, 
w)\cdot v\\
                    &  &\quad -u\cdot v\cdot w  - u\cdot w\cdot v  - 
w\cdot u\cdot v\   \\
                      &=& (u\cdot v) * w - u \cdot (v*w) - (u*w)\cdot v + 
u\cdot w \cdot v\ , \\
M_{12}(u, vw)&=& u*(v\cdot w)  - M_{11}(u, v)\cdot w -  v\cdot M_{11}(u, 
w)\\
                    &  &\quad  -u\cdot v\cdot w  -  v\cdot u\cdot w  -  
v\cdot w\cdot u\  \\
                      & =& u*(v\cdot w) - (u*v)\cdot w - v\cdot (u*w)+ 
v\cdot u\cdot w\ .\\
\end{eqnarray*}

\begin{prop}\label{infinit-2as}
The family of $(p+q)$-ary operations $M_{pq}$ constructed 
above defines  a  functor
$$(-)_{\Bi}: \{2as{\rm -algebras}\}\longrightarrow \{\Bi{\rm 
-algebras}\}.$$
\end{prop}
\P  Let us write $\Id_i$ for $\Id^{\t i}$. We want to prove that, for any integers $p,q,r\geq 1$, the 
following formula (denoted $\RRR_{pqr}$ in \ref{Binfini}) holds:

$$\sum _{(\uu p,\uu q)} M_{l r}\circ (M_{p_1 q_1}\cdots M_{p_l q_l} \t 
\Id_r)=\sum _{(\uu q,\uu r)} M_{p s}\circ (\Id_p
\t M_{q_1 r_1}\cdots M_{q_s r_s}),
$$
where the left sum is taken over all the partitions
$$(\uu p,\uu q)=(p_1,\cdots ,p_l;q_1,\cdots ,q_l)$$
  with $1\leq l\leq p+q$, $0\leq p_i\leq p$, $0\leq q_i\leq q$,
$\sum _ip_i=p$ and $\sum _iq_i=q$, and the right sum is taken over all 
partitions
$(\uu q,\uu r)=(q_1,\cdots ,q_s;r_1,\cdots ,r_s)$ with $1\leq s\leq 
q+r$, $0\leq q_j\leq q$, $0\leq r_j\leq r$,
$\sum _jq_j=q$ and $\sum _jr_j=r$.

We use the following notation: $\mu(x,y):= x*y,\  \nu(x,y):=x\cdot y$, 
$$\nu^{n-1}(u_1\ldots u_n) := u_1\cdot \ldots  \cdot u_n\  \textrm{ and }\ 
 \nu ^{p,q,r} := \nu^{p-1} \t \nu^{q-1}\t \nu^{r-1} \ .$$
We will prove, by induction on $h=p+q+r$ that the difference of the 
two sums is equal to
$$\mu \circ (\mu \t \Id_r) \circ ( \nu ^{p,q,r} ) - \mu \circ (\Id_p\t 
\mu ) \circ ( \nu ^{p,q,r} )$$
and so  is 0 by the associativity property of $\mu $.

Let us first check the property  for $p=q=r=1$, i.e.~$h=3$. Consider the 
two sums
\begin{eqnarray*}
S_l &:=& M_{21}(uv+vu, w) + M_{11}(M_{11}(u,v),w ) ,\\
S_r &:= &M_{12}(u, vw+wv) + M_{11}(u, M_{11}(v,w ))\ .
\end{eqnarray*}
By replacing the elements by their value in terms of the operations $*$ 
and $\cdot$ we get, up to a permutation of
variables, the following 8 different types of elements: 
$(-\circ_1-)\circ_2-$ and  $-\circ_1(-\circ_2-)$ for $\circ_i
= *$ or $\cdot$ . Collecting the terms of type $(-\cdot-)\cdot-$ (resp. 
  $-*(-\cdot-)$ ) in $S_l$ gives 0 and similarly in $S_r$.
Collecting the terms of type $(-*-)\cdot-$ (resp.  $-\cdot(-* -)$, 
resp.  $-\cdot (-*-)$, resp.  $-\cdot (-\cdot-)$ )  in $S_l$ gives the
same element as in  $S_r$. Collecting the terms of type $(-*-)*-$ or 
$-*(-*-)$ in $S_l$ gives $(u*v)*w$ and in $S_r$
gives $u*(v*w)$. Hence by the associativity property of $*$ we get 
$S_l=S_r$, i.e.~relation $\RRR_{111}$.

For $h\geq 4$, observe that the following equalities hold:
$$
\displaylines{
\mu \circ (\mu \t \Id_r) \circ  \nu ^{p,q,r} =\sum  _{(\uu p,\uu q)}
M_{l r}\circ (M_{p_1 q_1}\cdots M_{p_l q_l} \t \Id_r)+\hfill  \cr
  \nu \circ \big( \sum _{k+m+n<h}\ \mu \circ (\mu \t \Id_{r-m}) \circ \nu
^{p-k,q-m,r-n} \ \t \ \sum _{(\uu k,\uu m)} M_{l m}\circ
(M_{k_1 n_1}\cdots M_{k_l n_l} \t \Id_m)\big),\cr
}$$
and
$$
\displaylines{ \mu \circ (\Id_p\t \mu ) \circ ( \nu ^{p,q,r} )=\sum 
_{(\uu q,\uu r)} M_{p s}\circ (\Id_p \t M_{q_1 r_1}\cdots M_{q_s r_s})+\hfill \cr
  \nu \circ \big( \sum _{k+m+n<h} \mu \circ (\Id_{p-k}\t \mu ) \circ
\nu ^{p-k,q-m,r-n}
\t \sum _{(\uu n,\uu m)} M_{k s}\circ (\Id_k \t M_{n_1 m_1}\cdots
M_{n_s q_m})\big).\cr
}$$
The recursion hypothesis and the associativity of $\mu $ imply that
$$
\displaylines{
\nu \circ \big( \sum _{k+m+n<h}\ \mu \circ (\mu \t \Id_{r-m}) \circ \nu
^{p-k,q-m,r-n} \ \t \ \sum _{(\uu k,\uu m)} M_{l m}\circ
(M_{k_1 n_1}\cdots M_{k_l n_l} \t \Id_m)\big) \cr
\hfill =\nu \circ \big( \sum _{k+m+n<h} \mu \circ (\Id_{p-k}\t \mu ) \circ
\nu ^{p-k,q-m,r-n}
\t \sum _{(\uu n,\uu m)} M_{k s}\circ (\Id_k \t M_{n_1 m_1}\cdots
M_{n_s q_m})\big).\cr
}$$
So, taking the difference, we get
$$
\displaylines{
 \sum _{(\uu p,\uu q)} M_{l r}\circ (M_{p_1 q_1}\cdots M_{p_l q_l} \t 
\Id_r)-\sum _{(\uu q,\uu r)} M_{p s}\circ (\Id_p
\t M_{q_1 r_1}\cdots M_{q_s r_s})\hfill \cr
\hfill=\mu \circ (\mu \t \Id_r) \circ (\nu ^{p,q,r}) - \mu \circ (\Id_p\t \mu ) 
\circ (\nu ^{p,q,r}),\cr
}$$
as expected. \hfill $\square$

\begin{cor}\label{cor:morphism}
Let  $R$ be a $\Bi$-algebra and let 
$(T^c(R), *, \DD)$ be the associated bialgebra. Put on $T^c(R)$ the $\Bi$-algebra  structure induced by the 2-associative operations $*$ and $\cdot$ (concatenation). Then the inclusion $R\to 
T^c(R)$ is a $\Bi$-morphism. In other words the $\Bi$-algebra structure of $R$ has been extended to $T^c(R)$.
\end{cor}

\P The operation $*$ is given by Proposition \ref{Tfc} and the operation $\cdot$ is concatenation. Denote 
 by  $M_{pq} $ the
  $\Bi$-structure on $R$ and by ${\overline M}_{pq} $ the $\Bi$-structure on $T^c(R)$ given by Proposition \ref{infinit-2as}. Let us show that ${\overline M}_{pq} =   M_{pq} $ on $R^{\t p+q}$ by induction on $p+q$. We use
formula (\ref{Mpqinductif}) which holds for ${\overline M}$ and $M$ on $R$. By 
induction ${\overline M} = M $ on the right side, therefore
${\overline M}_{pq} = M_{pq} $ for any $(p,q)$. \hfill $\square$

\subsection{Universal enveloping \dda algebra}\label{uea}
 We define the {\it \ue \dda algebra} of a $\Bi$-algebra $R$, denoted
  $U2(R)$, as the
quotient of the free \dda algebra $2as(R)$ over the vector space $R$ by 
the relations
$$M_{pq}(r_1\cdots r_p , s_1\cdots s_q) \approx \widetilde 
M_{pq}(r_1\cdots r_p , s_1\cdots s_q),
\  r_i, s_j \in R$$
where $M_{pq}$ denotes the operation in $R$ and  $\widetilde M_{pq}$ 
denotes the operation in $2as(R)$:
$$U2(R):= 2as(R)/ \approx\ . $$
In other words we divide $2as(R)$ by the ideal (in the 2-associative algebra sense) generated by the elements
$$M_{pq}(r_1\cdots r_p , s_1\cdots s_q) - \widetilde M_{pq}(r_1\cdots 
r_p , s_1\cdots s_q) ,$$
  thus $U2(R)$ is a 2-associative algebra.

\begin{lemma}\label{lemma:left adjoint}
The functor $U2:  \{\Bi{\rm 
-alg}\}\longrightarrow \{2as{\rm -alg}\} $ is left
adjoint to the functor $(-)_{\Bi}: \{2as{\rm -alg}\}\longrightarrow 
\{\Bi{\rm -alg}\}$ .
\end{lemma}

\P Let $A$ be a $2as$-algebra and let $f: R\to A_{\Bi}$ be a
morphism of $\Bi$-algebras. It determines uniquely a morphism of \dda 
algebras $2as(R) \to A$ since $A$ is a
\dda algebra and $2as(R)$ is free. This morphism passes to the quotient 
to give $U2(R) \to A$ because the image in $A$ of the two elements 
\[M_{pq}(r_1\cdots r_p , s_1\cdots s_q) \textrm{ and } \widetilde M_{pq}(r_1\cdots 
r_p , s_1\cdots s_q) ,\]
 is the same, namely 
$M_{pq}(f(r_1)\cdots f(r_p) , f(s_1)\cdots f(s_q))$.

On the other hand, let $g: U2(R) \to A$ be a morphism  of \dda 
algebras. From the construction of $U2(R)$ it follows that the map $R \to U2(R)$ is a $\Bi$-morphism. Hence the composition with $g$ gives  a $\Bi$-morphism $R\to A$.

Clearly these two constructions are inverse to each other, and 
therefore $U2$ is left adjoint to $(-)_{\Bi}$. \hfill $\square$

\begin{cor}\label{cor:free}
 The \ue \dda algebra of the free 
$\Bi$-algebra is canonically isomorphic to the free \dda algebra:
\[ U2(\Bi(V))\cong 2as(V)\ .\]
\end{cor}
\P First recall that $A_{\Bi}$ has the same underlying vector space as 
$A$. Since $U2$ is left adjoint to
$(-)_{\Bi}$ and since $\Bi$ is left adjoint to the forgetful functor, 
the composite is left adjoint to the forgetful
functor from \dda algebras  to vector spaces. Hence it is the functor 
$2as$.
\hfill $\square$

\subsection{Definition}\label{def:dda bialgebra}
{\rm A {\it  \dda  bialgebra} (resp.~{\it \dda 
Hopf algebra})  $(\HH, *,\cdot, \DD)$ is a vector space $\HH$
equipped with two operations $*$ and $\cdot$ and one cooperation $\DD$, 
such that

\noindent $\bullet$ $(\HH, *,\DD)$ is a bialgebra (resp.  Hopf algebra), 
cf.~\ref{Hopf},

\noindent $\bullet$ $(\HH,\cdot, \DD)$ is a unital infinitesimal bialgebra, 
cf.~\ref{Def:infinitesimal}.}

\begin{prop}\label{prop:primitive}
For any \dda bialgebra (e.g. ~\dda Hopf 
algebra) its primitive part is a
sub-$\Bi$-algebra.
\end{prop}

\P
Given elements $x_1,\dots ,x_n,y_1,\dots ,y_m$ in $\Prim \HH $, we need 
to prove that $M_{n m}(x_1\dots x_n,y_1
\dots y_m)$ belongs to $\Prim \HH $ too, that is $\DDo \circ M_{nm}=0$ 
on $\Prim \HH$. The proof is by induction on $(n,m)$.
Instead of giving the details we give an alternate proof in the case 
where $\HH$ is connected. We only need this case in the
sequel of the article.

\noindent If $\HH$ is a connected 2-associative bialgebra, then $\HH$ is 
isomorphic to $T^{fc}(\Prim \HH)$ by Theorem \ref{thm:infinitesimal} and $\Prim \HH$ is a
$\Bi$-algebra. So we can apply Corollary \ref{cor:morphism}. \hfill $\square$

\begin{prop}\label{prop:unique}
There exists a unique cooperation $\DD$ 
on the free \dda algebra $2as(V)$ which makes it into a
\dda bialgebra and for which $V$ is primitive. As a coalgebra $2as(V)$ is connected.
\end{prop}

\P We define $\DD : 2as(V) \to 2as(V)\t 2as(V) $ by the following 
requirements:
\\

\noindent $\bullet\quad  \DD(1) = 1\t 1 ,$

\noindent $\bullet\quad \DD(v) = v\t 1 + 1\t v ,$ for $v\in V$,

\noindent $\bullet\quad \DD(x*y) = \DD(x)*\DD(y)  ,$

\noindent $\bullet\quad \DD(x\cdot y) = (x\t 1)\cdot \DD(y) + \DD(x)\cdot (1\t 
y) - x\t y  .$
\\

Let us prove that $\DD$ is well-defined by induction on the degree of 
the elements in
$2as(V)=\oplus_{n\geq 0}2as(V)_n$. It is already defined on
$2as(V)_0= K.1$ and $2as(V)_1=V$. Suppose that $\DD$ is defined up to  
$2as(V)_{n-1}$. We define it on
$2as(V)_n$ as follows. Any element of $2as(V)_n$ is of the form $x*y$ 
or $x\cdot y$ for elements $x$ and $y$ of
degree strictly smaller than $n$. Then
$\DD(x*y) $ and $ \DD(x\cdot y) $ are given by the required formulas. 
Since the only relations are the associativity of $*$, the
associativity of $\cdot$ and the unitality of 1 for both products, we 
need to verify that
\begin{eqnarray*}
  \DD((x* y) *z) =  \DD(x*(y * z)) , \\
  \DD((x\cdot y) \cdot z) =  \DD(x\cdot (y \cdot z)) , \\
  \DD(1*x) =  \DD(x) =  \DD(x*1) , \\
  \DD(1\cdot x) =  \DD(x) =  \DD(x\cdot 1) . 
\end{eqnarray*}
The last two lines are straightforward to check. The first one is 
classical. Let us check the second one. On the left side
we get (with the notation $a\cdot \DD(b) = (a\t 1) \cdot \DD(b)$ and 
$\DD(a)\cdot b = \DD(a)\cdot (1\t b)$ )
\begin{eqnarray*}
\DD((x\cdot y) \cdot z) &=& (x\cdot y) \cdot  \DD(z) + \DD(x\cdot 
y)\cdot z - (x \cdot y)\t z ,\\
&=& x\cdot y\cdot  \DD(z) + (x\cdot \DD(y)+ \DD(x)\cdot y - x\t y)\cdot 
z - x \cdot y\t z ,\\
&=& x\cdot y\cdot  \DD(z) + x\cdot \DD(y)\cdot z+ \DD(x)\cdot y\cdot z - 
x\t y\cdot z - x \cdot y\t z \ .
\end{eqnarray*}
On the right side we get
\begin{eqnarray*}
\DD(x\cdot (y\cdot z))&=&x\cdot \DD(y\cdot z) + \DD(x)\cdot (y\cdot z) - 
x\t y\cdot z ,\\
&=& x\cdot (y\cdot \DD(z) + \DD(y)\cdot z - y \t z)+ \DD(x)\cdot y\cdot 
z - x\t y\cdot z  ,\\
&=& x\cdot y\cdot  \DD(z) + x\cdot \DD(y)\cdot z+ \DD(x)\cdot y\cdot z - 
x\t y\cdot z - x \cdot y\t z \ .
\end{eqnarray*}
We see that they are equal. Hence $\DD$ is well-defined.
There is a unique map $2as(V) \to 2as(V)\t 2as(V)\t 2as(V)$ sending $1$ to $1\t 1\t 1$ and $v$ to $v\t 1\t 1 + 1 \t v\t 1 + 1 \t 1 \t v$ 
and compatible with $*$ and $\cdot$. Since both maps $(\DD\t \Id)\circ \DD$ and $(\Id \t \DD) \circ \DD$ satisfy these properties, they 
are equal, and so $\DD$ is coassociative.

 It follows  from this computation that $(2as(V), *, \cdot , \DD)$ is a \dda bialgebra.

\noindent Proof of connectedness. It is sufficient to prove that any element 
in $2as(V)/K$ is conilpotent. An element in $V$
is primitive, hence conilpotent. So it suffices to prove that, if $x$ 
and $y$ are conilpotent, then so are $x*y$ and $x\cdot
y$. This fact follows from the relationship between $\DD$ and $*$, 
respectively $\DD$ and $\cdot$.\hfill $\square$

\begin{cor}\label{cor:bialgebra}
For any $\Bi$-algebra $R$, $U2(R)$ is a 
connected 2-asso-\allowbreak ciative bialgebra.  \hfill $\square$
\end{cor}

\section{The main theorem}\label{S:main}
In this section we state and prove the main results, that is the 
structure theorem for cofree Hopf algebras and the
main theorem from which it follows, that is the structure theorem for 
connected \dda bialgebras. They are analogue of
the Cartier-Milnor-Moore theorem (for (a)$\Rightarrow$ (b)) and the Poincar\'e-Birkhoff-Witt theorem
(for (b)$\Rightarrow$ (c)),  which hold for {\it cocommutative} connected
bialgebras (over a characteristic zero field). Observe that in the 
non-cocommutative setting we do not need the characteristic zero 
hypothesis.

\subsection{Classification of cofree bialgebras}\label{classification} A connected cofree bialgebra can be equipped with a second product 
(by using cofreeness) which makes it into a connected 2-associative bialgebra. Our main result about these objects is the following.

\begin{thm}\label{thm:main} If $\HH$ is a bialgebra over a field $K$, 
then the following are equivalent:

\noindent (a) $\HH$ is a connected \dda  bialgebra,

\noindent (b) $\HH$ is isomorphic to $U2(\Prim \HH)$ as a \dda bialgebra,

\noindent (c) $\HH$ is cofree among the connected coalgebras.
\end{thm}

Since a connected bialgebra is a Hopf algebra, one can 
replace connected bialgebra by
connected Hopf algebra in this theorem.
\\

\P  We prove the following implications (a) $\Rightarrow$ (c) 
$\Rightarrow$ (b) $\Rightarrow$ (a). We put $R:=\Prim \HH$.

\noindent (a) $\Rightarrow$ (c). If $\HH$ is a connected \dda  bialgebra, then, 
by Theorem \ref{thm:infinitesimal}, $\HH$ is isomorphic to $T^{fc}(R)$ as
a unital infinitesimal bialgebra. Therefore $\HH$ is cofree.

\noindent (c) $\Rightarrow$ (b).  If $\HH$ is cofree, then it is isomorphic to 
$T^{fc}(R)$ and $R$ is a $\Bi$-algebra. Observe that $T^{fc}(R)$ is a $2as$-algebra: the product $*$ is inherated from the 
associative product of $\HH$ under the isomorphism and the product $\cdot$ is the concatenation. By adjunction 
(cf.~Lemma \ref{lemma:left adjoint}), the inclusion map $R\to T^{fc}(R)$ gives rise to a $2as$-morphism
$$ \phi : U2(R) \to T^{fc}(R)\ .$$
On the other hand, the inclusion $R \to U2(R)$ admits a unique extension
$$ \psi : T^{fc}(R)\to U2(R)$$
by $\psi(r_1 \ldots r_n ) = \psi(r_1) \cdot \ldots \cdot 
\psi(r_n)$. It is immediate to check that $\phi$ and
$\psi$ are inverse to each other.

\noindent (b) $\Rightarrow$ (a). This is Corollary \ref{cor:bialgebra}.
\hfill $\square$

As an immediate consequence, we obtain a structure 
theorem for cofree bialgebras:

\begin{cor} There is an equivalence between the 
category of cofree bialgebras
and the category of bialgebras of the form $U2(R)$ for some 
$\Bi$-algebra $R$.
\end{cor} \hfill $\square$

The particular case of free algebras reads as follows:

\begin{cor}\label{isofree} For any vector space $V$  the free \dda 
algebra $2as(V)$ and the free $\Bi$-algebra $\Bi (V)$ are
related by the following isomorphisms:
$$\bar \phi :\Prim (2as(V)) \cong \Bi (V)\ \hbox{as $\Bi$-algebras},$$
$$\phi : 2as(V) \cong   T^{fc}(\Bi (V))\  \hbox{as \dda bialgebras}.$$
\end{cor} 
{}\hfill $\square$

\subsection{The  differential graded framework}\label{differential}
So far we worked in the category of vector spaces over $K$. However we can, more generally, work in the category of graded vector spaces, or differential graded vector spaces (that is chain complexes), in which all the constructions and results remain true. In the graded case the formulas are exactly the same as in the non-graded case, provided we write them as equalities among maps, that is we forget the entries. For instance the relation  $(\RRR_{111})$ should be written:
$$M_{21}\circ ( (12)+1_3)+ M_{11}\circ (M_{11}\t 1_1) =  M_{12}\circ (1_3+(23))+M_{11}\circ (1_1\t M_{11})
\ .$$
In this formula $1_3$ is the identity on $V^{\t 3}$ and $(12)$ is the map which permutes the first two variables, so
$$(12)(uvw) = (-1)^{\vert u\vert \vert v\vert} vuw$$
where $\vert u\vert$ is the degree of the homogeneous element $u$. One should also use the Koszul sign rule, that is $f\t g$ means the map which sends $u\t v$ to $ (-1)^{\vert g\vert \vert u\vert} f(u)\t g(v)$.

In the differential graded case, a $\Bi$-algebra structure on the graded vector space $V=\oplus_{n\in {\mathbf Z}} V^{\t n}$ is equivalent to a graded differential algebra structure on the tensor coalgebra $T(V[1])$ where $V[1]$ is the desuspension of $V$ (cf. ~\cite{GJ} or \cite{Vor} for details).
\\

\subsection{The dual  framework}\label{dual}

Instead of starting with a bialgebra which is cofree as a coalgebra, we could start with a bialgebra which is free as an algebra. Then similar results hold, but 2-associative algebras have to be replaced by 2-coassociative coalgebras and so forth. Observe that the notion of unital infinitesimal bialgebra is self-dual (this is easily seen on the picture of \ref{Def:infinitesimal}). See \ref{selfduality} for the free case.
\\
\newpage

\section{Free \dda algebra}\label{S:free2as}
In this section we give an explicit description of the free \dda algebra in terms of planar trees. The free $2as$-algebra on one generator can be identified with the non-commutative polynomials over the planar (rooted) trees. Alternatively the underlying vector space can be identified with the vector space spanned by 
{\it two} copies of the set of planar trees amalgamated over the 
trivial tree.

We study the operad  associated to \dda algebras, we compute its dual (in the operadic sense) 
and we show that they are Koszul operads. As a
consequence we obtain a chain complex to compute the (co)homology of 
\dda algebras. 
\subsection{Planar trees}\label{planartree}
We denote by $\TT_n$ the set of {\it planar (rooted)
trees} with $n$ leaves,
$n\geq 1$ (and one root) such that each internal vertex has one root and at least two offsprings. We say that such a vertex has valence at least 2.
Here are the first sets $\TT_n$:
$$\TT_1 = \{\vert\} ,\qquad \ \TT_2 = \{ \arbreA  \} ,\qquad \TT_3 = 
\{ \arbreAB , \arbreBA ; \arbreBB \}$$
%$$\TT_4 = \{\ \arbreun, \arbredeux, \arbretrois, \arbrequatre, 
%\arbrecinq;\arbreutt , \arbreuut, \arbretut, \arbretuu, \arbrettu; \arbrettt \}. $$
The lowest vertex of a tree is called the {\it root vertex}. An edge which is 
neither the root nor a leaf is called an {\it internal
edge}. Observe that the trivial tree $\mid$ has no vertex.

The integer $n$ is called the {\it degree} of $t\in \TT_n$.
We define $\TT_{\infty} := \bigcup_{n\geq 1} \TT_n$. The number of 
elements in $\TT_n$ is the so-called \emph{super
Catalan number} or \emph{Schr\" oder number}, denoted 
$C_{n-1}$ whose value is :
$$1, 1, 3, 11, 45, 197, 903,  \ldots\ .$$
It is well-known that its generating series is
\begin{equation}\sum_{n\geq 0} C_nx^{n}= {(1+x - \sqrt {1-6x+x^2})}/{4x} 
\end{equation}
(see also \ref{generating}).

By definition the {\it grafting} of $k$ planar trees $\{t^{1}, \ldots , 
t^{k}\}$ is a planar tree denoted $t^{1}\vee
\ldots
\vee t^{k}$ obtained by joining the $k$ roots to a new vertex and 
adding a new root.
Any planar tree can be
uniquely obtained as
$t=t^{1}\vee \ldots \vee t^{k}$, where $k$ is the valence of the root 
vertex (for $t=\mid$, one has $k=1$ and
$t^1=\mid$). Observe that the grafting operation is not associative.
At some point in the sequel we will need to work with two copies of 
$\TT_{\infty}$ that we denote by
$\TT_{\infty}^{*}$ and $\TT_{\infty}^{\bullet}$. Pictorially we identify 
their elements by decorating the root vertex by
$*$ or
$\cdot$ respectively. By convention we identify the copy of the tree 
$\mid$ (which has no root vertex) in
$\TT_{\infty}^{*}$ to its avatar in $\TT_{\infty}^{\bullet}$.

\subsection{Free \dda algebra}\label{free2ass} 
Since, in the relations defining the 
notion of \dda algebra, the variables stay in the same
order, we only need to understand the free \dda algebra on one 
generator (i.e.~over $K$). In operadic terminology we are
dealing with a regular operad (i.e.~constructed out of a non-$\SS$-operad).

More explicitly, if we denote by $2as(K)= \oplus_{n\geq 0}2as_n$ the 
free 2-associative algebra on one generator, then the
free 2-associative algebra on $V$ is
$$2as(V) =  \bigoplus_{n\geq 0}2as_n\t V^{\t n}\ .$$
The \dda structure of $2as(V)$ is induced by the \dda structure of 
$2as(K)$ and concatenation of tensors:
\begin{eqnarray*}
(s; v_1\cdots v_p)*(t, v_{p+1}\cdots v_n) & = & (s*t; v_1\cdots v_n),\\
(s; v_1\cdots v_p)\cdot (t, v_{p+1}\cdots v_n) & = & (s\cdot t; v_1\cdots v_n).
\end{eqnarray*}

Let $T(\TTi) = (T(K[\TTi]), *)$ be the free unital associative algebra over the vector 
space $K[\TTi]$ generated by the set
$\TTi$. So here $*$ is the concatenation. A set of linear  generators of $T(\TTi)$ is made of the (non-commutative) monomials
$t_1*t_2*\cdots *t_k$ where the $t_i$'s are planar trees. We
define the product $\cdot $ as follows:
\begin{eqnarray*}
s \cdot  t &:=& s^1 \vee \ldots \vee s^m\vee t^1\vee \ldots
\vee t^n ,\\
(s_1* \ldots *s_k)\cdot  t &:=& ( s_1 \vee \ldots \vee s_k)\vee t^1\vee 
\ldots
\vee t^n, \quad \hbox{ for } k\geq 2, \\
s\cdot  (t_1*\ldots *t_l) &:=& s^1 \vee \ldots \vee s^m\vee (t_1\vee 
\ldots
\vee t_l), \quad \hbox{ for } l\geq 2 ,\\
(s_1* \ldots *s_k)\cdot  (t_1*\ldots *t_l) &:=& ( s_1 \vee \ldots \vee 
s_k)\vee (t_1\vee \ldots
\vee t_l), \quad \hbox{ for } k\geq 2, l\geq 2 ,
\end{eqnarray*}
where $s=s^1\vee \ldots \vee s^m$ and $t= t^1\vee \ldots \vee t^n$.

Observe that the dot product of two monomials is always a tree.

\begin{lemma}\label{lemma:bitree} There is a natural bijection $T(\TTi)\to 
K 1\oplus K[\TTi^{*}\bigcup \TTi^{\bullet}]$ obtained by
sending a tree in $\TTi$ to its counterpart in
$\TTi^{\bullet}$ and by sending a nontrivial monomial $t_1*\cdots *t_n$ 
(i.e.~$n\geq 2$) to the grafting $t_1\vee\ldots
\vee t_n$ in $\TTi^{*}$. Hence we have $\dim T(\TTi)_n = 2C_{n-1}$ for 
$n\geq 2$.
\end{lemma}

Explicitly we obtain:
\begin{displaymath}\begin{array}{cc}
s^* * t^* = (s_1\vee \ldots \vee s_n\vee t_1\vee \ldots \vee t_m)^*, & 
s^* \cdot t^* =
(s\vee  t)^{\cdot},\\
s^{\cdot} * t^* = (s\vee t_1\vee \ldots \vee t_m)^*, & s^{\cdot} \cdot 
t^* =
(s_1\vee \ldots \vee s_n\vee t)^{\cdot},\\
s^* * t^{\cdot} = (s_1\vee \ldots \vee s_n\vee t)^*, & s^* \cdot 
t^{\cdot} =
(s\vee t_1\vee \ldots \vee t_m)^{\cdot},\\
s^{\cdot} * t^{\cdot} = (s\vee t)^*,& s^{\cdot} \cdot t^{\cdot} = 
(s_1\vee \ldots
\vee s_n\vee t_1\vee \ldots \vee t_m)^{\cdot},
\end{array}\end{displaymath}
whenever $s=s_1\vee \ldots \vee s_n$ and $ t= t_1\vee \ldots \vee t_m$. 
\\

\P The grafting operation gives a bijection between the set of $n$-tuples 
of planar trees, for any $n\geq 2$, and
$\TT_{\infty}\backslash \{\mid\}$. Indeed the inverse is the map 
$t_{1}\vee \cdots\vee
t_{n}\mapsto\{t_{1}, \cdots , t_{n}\}$. It  gives a bijection between the monomials of 
degree $\geq 2$ in $\TTi$ and $K[\TTi^{\bullet}]$. \hfill $\square$

\begin{thm}\label{thm:free2as} The free \dda algebra $2as(K)$ on one 
generator is isomorphic to the
\dda algebra $(T(\TTi), *, \cdot)$ where $*$ is the concatenation product in $T(\TTi)$ and $\cdot$ is given by the formulas in \ref{free2ass}.
\end{thm}

\P First we check that the dot product is associative. There are 
eight cases, which are immediate to check. For
instance, if
$u = u^1\vee \ldots \vee u^p$, then
$$(s\cdot t)\cdot u =  s^1 \vee \ldots \vee s^m\vee t^1\vee \ldots
\vee t^n  \vee u^1\vee \ldots
\vee u^p =s\cdot (t\cdot u)\ .$$

Second we prove that, as a $2as$-algebra, it is free on one generator. Let $x$ be the generator of the free \dda algebra $2as(K)=  
\oplus_{n\geq 0} 2as_n$. Since $T(\TTi)$ is a
\dda algebra, mapping $x$ to the tree $\mid$ induces a \dda 
morphism $\aa : 2as(K) \to T(\TTi)$.
Since any linear generator of
$T(\TTi)$ can be obtained from $\mid$ by using $*$ and $\cdot$, this 
map is surjective. The degree $n$ part of
$T(\TTi)$ is of dimension $2C_{n-1}$ by Lemma \ref{lemma:bitree}, hence to prove that 
$\aa$ is an isomorphism it suffices to prove that
the dimension of $2as_n$ is less than or equal to $2C_{n-1}$.

We identify a parenthesizing of a word to a planar \emph{binary} 
tree. The space $2as_n$ is linearly generated by the planar binary
trees with $n$ leaves such that each vertex is decorated
by either $*$ or $\cdot$ . Because of the associativity relations of 
both operations, two trees which differ  only
by replacing locally
$$\arbreAB \qquad \hbox {by }\qquad  \arbreBA$$
where the vertices are all decorated by $*$, or all decorated by $\cdot$ ,
give rise to the same element in $2as_n$. Let us draw it as a tree 
obtained from the previous one by replacing the
equivalent patterns by

$$\arbreBBe  \qquad \txt{or}\qquad \arbreBBdot \qquad \txt{respectively.}$$
Applying the process ad libitum gives a planar tree such that each 
internal edge has different decorations at both ends.
So, the decorations are completely determined by the decoration of the 
root-vertex, and so we can ignore the other ones.
This argument shows that the dimension of $2as_n$ is less than or equal 
to $2C_{n-1}$ as expected.

  So, in fact, the dimension of
$2as_n$ is precisely $2C_{n-1}$ for $n\geq 2$. \hfill $\square$

\begin{cor}\label{cor:free2as} The following is the free \dda algebra on one generator (namely $\vert$ ). As a vector space it is $K 1\oplus 
K[\TTi^{*}\bigcup \TTi^{\bullet}]$ where the two copies of $\vert$ have been identified. The laws are given by
\[ s^{\circ_s} \circ t^{\circ_t} = ((s\vee t)/\sim )^{ \circ}\]
where ${\circ_s}, \circ_t$ and $ \circ$ are either $*$ or $\cdot$, and the quotient $\sim$ is as follows: the edge from the root of $s$ (resp.~$t$) to the new root is collapsed to a point if $\circ_s =  \circ$ (resp. $\circ_t =  \circ$). 
\end{cor} \hfill  $\square$

\subsection{Remark} \label{rem:freeness} {\rm One can switch the roles of the two products, in particular 
$2as(V)$ is free as an associative algebra for $\cdot$ .

A free 2-associative set is called a duplex in  \cite{Pir}.
Corollary \ref{cor:free2as} gives also the structure of the free duplex in one generator.

Observe that for $s=s_1\vee \ldots \vee s_n$ one has 
$$s^* = s_1^{\cdot}*\ldots * s_n^{\cdot}\quad \hbox{ and }\quad 
s^{\cdot}=  s_1^*\cdot\ldots \cdot s_n^*\ .$$
In figure 1 we draw the elements in low dimension.

\begin{figure}
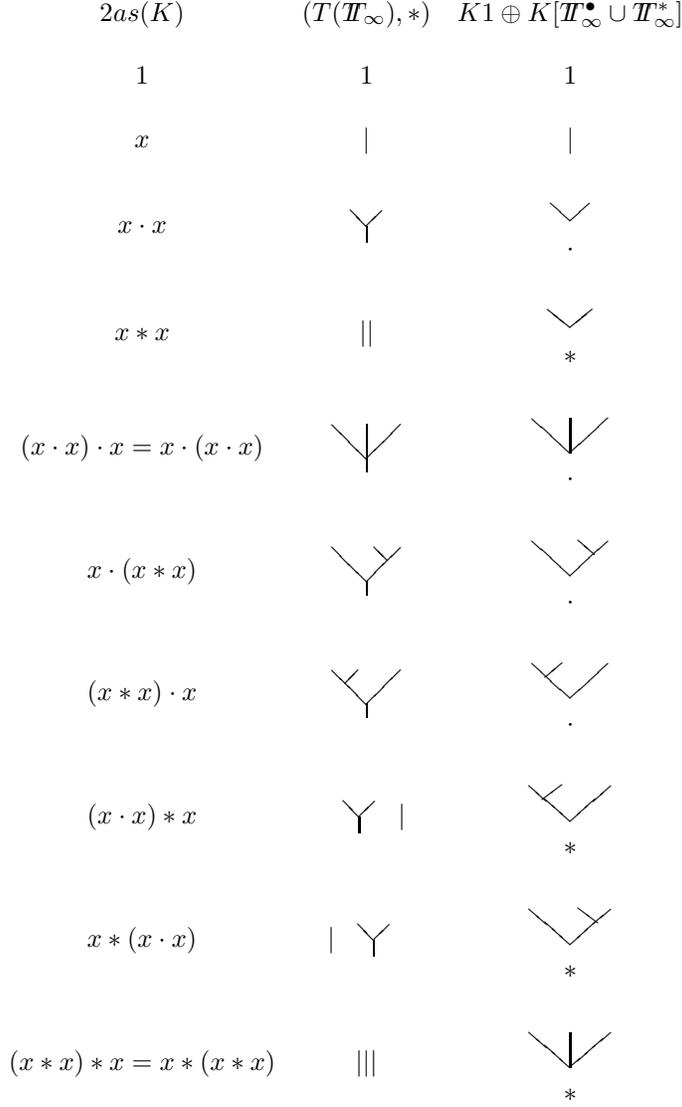
\begin{displaymath}\begin{array}{ccc}
2as(K)   & (T( \TTi ),*) & K1\oplus K[\TTi^{\bullet} \cup \TTi^{*} ]\\
&&\\
1 & 1 & 1  \\
&&\\
x & \mid & \mid  \\
&&\\
x\cdot x &\arbreA & \arbreAdot \\
&&\\
x*x &\mid \mid & \arbreAe \\
&&\\
(x\cdot x)\cdot x = x\cdot (x\cdot x)&\arbreBB &  \arbreBBdot \\
&&\\
x\cdot (x* x) &\arbreAB &  \arbreABdot\\
&&\\
(x* x)\cdot x  &\arbreBA &  \arbreBAdot\\
&&\\
(x\cdot x)* x &\arbreA\mid  &  \arbreABe\\
&&\\
x*(x\cdot  x) &\mid \arbreA &  \arbreBAe\\
&&\\
(x*x)*x = x*(x*x) & \mid \mid \mid &  \arbreBBe 
\end{array}\end{displaymath}\label{figure}\caption{Low dimensional elements in $2as(K)$}
\end{figure}

}

\subsection{Involution}\label{involution}
Let us define a map $\iota : 2as(V) \to 2as(V)$ by the following requirements:
$$\iota (1) = 1,\ \iota (v)=v,\ \iota (x*y) = \iota (y)*\iota (x),\ \iota (x\cdot y) = \iota (y)\cdot \iota (x).$$
First, it is immediate to verify that $\iota ^2 = \Id$. Second, by using the compatibility relations between $\DD$ and the products one can show that
$$\DD(\iota (x)) = \iota (x_{(2)})\t \iota (x_{(1)}), {\rm where}\ \DD(x) = x_{(1)}\t x_{(2)}. $$
On the trees this involution is simply the symmetry around the root axis.

\subsection{Generating series}\label{generating} Let $C(x):=\sum_{n\geq 0} C_nx^{n}$ be the 
generating series for the super Catalan numbers. For a
graded vector space $V = \oplus_{n\geq 0} V_n$ we denote by 
$$f^V(x):= 
\sum_{n\geq 0} \dim V_n\, x^{n}$$
 its generating
series.  Therefore the generating series of $K[\TTi]$ is $f^{K[\TTi]}(x)= x 
C(x)$, while the generating series of $2as(K)$ is
$f^{2as}(x) = 1-x+2x C(x)$. Since $2as(K)=T(K[\TTi])$ one has 
$f^{2as}(x) = (1-f^{K[\TTi]}(x) )^{-1}$. It follows that $C(x)$
satisfies the algebraic relation $2xC(x)^2 - (1+x) C(x) + 1 = 0$. We 
deduce from it the expression given in \ref{planartree}:
$$\sum_{n\geq 0} C_nx^{n}= {(1+x - \sqrt {1-6x+x^2})}/{4x}\ .$$

\subsection{Homology and Koszul duality} In \cite{GK} Ginzburg and Kapranov developed the theory of Koszul duality for binary quadratic operads. When an operad is Koszul, its Koszul dual permits us to construct a small chain complex to compute the homology of algebras and also to construct the minimal model, hence giving rise to the algebras up to homotopy. This theory is applicable here and we describe  the outcome.

Until the end of this section we work in the nonunital framework, that is we do not suppose that the \dda algebras have a unit. In particular the degree zero component of the free \dda algebra is 0. Let us recall that the Hochschild chain complex of a nonunital associative algebra $A$ is of the form
$$ \cdots \to A^{\t n} {\buildrel {b'}\over {\longrightarrow}} A^{\t n-1}\to \cdots \to A$$
where $b'$ is given by $b'(a_1,\cdots , a_n) = \sum _{i=1}^{i=n-1}(-1)^i (a_1,\cdots , a_ia_{i-1},\cdots , a_n)$. Its homology is denoted $H^{As}_*(A)$.

From the definition of the dual of an operad (cf.~loc.cit. or \cite{Lod1} Appendix B, for a brief introduction), the dual of the operad $2as$ is the operad $2as^!$ whose algebras have two associative operations $*$ and $\cdot$ verifying:
\begin{eqnarray*}
(x*y)\cdot z&=\ 0\ =&x*(y\cdot z)\\
(x\cdot y)* z&=\ 0\ =&x\cdot (y*z)
\end{eqnarray*}
It is clear that the free $2as^!$-algebra over $V$ is $V$ in dimension 1 and the sum of two copies of $V^{\t n}$ in dimension $n\geq 2$, the first one corresponding to $*$ and the second one to $\cdot$ . Hence the chain complex of a \dda algebra is two copies of the Hochschild complex described above, amalgamated in dimension 1. To check the Koszulity of the operad $2as$ it suffices to show that the chain complex of the free $2as$-algebra is acyclic. Since $2as(V)$ is free as an associative algebra for $*$ by Theorem \ref{thm:free2as},  and also free as an associative algebra for $\cdot$ (cf.~\ref{rem:freeness}), the corresponding Hochschild complexes are acyclic and hence we have proved the following:

\begin{prop} The operad of \dda algebras is a Koszul operad.
\end{prop}\hfill $\square$

As a consequence, for a $2as$-algebra $(A,*,\cdot)$  we have, for $n\geq 3$, 
$$H^{2as}_n(A) = H^{As}_n(A,*)\oplus H^{As}_n(A,\cdot)\ .$$

\section{Bialgebra structure of the free $2as$-algebra and free $\Bi$-algebra}\label{S:6}

We unravel the Hopf algebra structure of the free \dda algebra. 
 From Theorem \ref{thm:main} and the explicitation of the free 2-associative algebra we deduce the structure of the operad $\Bi$. 

\subsection{The coproduct on $2as(V)$}\label{coproduct2as}  Proposition \ref{prop:unique} determines a unique coproduct $\DD$ on $2as(V)$, for which\\
 $(2as(V), *, \DD)$ is a Hopf algebra. The following Proposition gives a recursive formula for $\DD$ in terms of $V$-decorated trees. By definition a \emph {$V$-decorated} tree is an element $(t; v_1 \ldots v_n)\in {\TT}_n\times V^{\t n}$. It is helpful to think of the entry $v_i$ as a decoration of the $i$th leaf of the tree $t$. 
\begin {prop}\label{prop:Delta} 
The coproduct $\Delta$ on the free $2$-associative algebra $T(\TTi)$ is given recursively by the following formula
$$\DD(t^1\vee\ldots \vee t^r) =
 1\t t + \sum_{i=1}^{r}\bigvee (t^1, \ldots , t^{i-1})\cdot {t^i}_{\{1\}} \t  {t^i}_{\{2\}}\cdot\bigvee (t^{i+1},\ldots , t^r) ,$$
where 
\begin{displaymath}
\bigvee (t^1, \ldots , t^{k}):= \left\{ 
\begin{array}{ll}
t^1\vee  \ldots \vee t^{k}& \textrm{if $ k>1$,} \\
t^{1(1)}  \ldots   t^{1(n)}& \textrm{if $k=1$ and $t^1 = t^{1(1)}\vee   \ldots \vee  t^{1(n)}$,}
\end{array}\right.
\end{displaymath} 
and we use the notation $s_{\{1\}}\t s_{\{2\}} := \DD(s^1\cdots s^k) - 1\t s^1\cdots s^k $ for $s = s^1\vee \cdots\vee  s^k$, and the convention $1\vee \oo = \oo$.
\end{prop}

\P It is straightforward (though tedious) to check by induction that this map $\DD$ satisfies the four identities explicitly written at the beginning of the proof of Proposition \ref{prop:unique}. For instance we have
$$\DD(\vert; v) = 1\t (\vert; v) + (\vert; v) \t 1.$$
{}\hfill $\square$

\subsection{Example} 
\[\begin{array}{l}
\bar \DD (\arbreA ; uv) = (\vert; u) \t (\vert; v) ,\\
\bar \DD (\arbreBA ; uvw) = (\vert; u)(\vert; v) \t (\vert; w) + (\vert; u)\t (\arbreA ; vw)  + (\vert; v)\t (\arbreA ; uw) , \\
\bar \DD (\arbreABC ; uvwx) = (\arbreA; uv)(\vert ; w) \t (\vert; x) +(\arbreA ; uv)\t (\arbreA; wx)\\
\hfill +(\vert ; w)\t (\arbreBB; uvx)
 + (\vert; u) \t (\arbreBA; vwx) + (\vert;u)(\vert,w)\t(\arbreA;vx).
\end{array}\]

\subsection{Selfduality and rigidity of $2as(V)$}\label{selfduality}
On the \dda bialgebra $(2as(V), *, \cdot, \DD)$ one can construct another coassociative cooperation $\dd$ by the following requirements:
\\

\noindent $\bullet\quad  \dd(1) = 1\t 1 ,$

\noindent $\bullet\quad \dd(v) = v\t 1 + 1\t v ,$ for $v\in V$,

\noindent $\bullet\quad \dd(x*y) =  (x\t 1)*\dd(y) + \dd(x)* (1\t 
y) - x\t y  ,$

\noindent $\bullet\quad \dd(x\cdot y) = \dd(x)\cdot\dd(y) .$
\\

In other words, the compatibility relation between the products $*, \cdot$ and the coproducts $\DD, \dd$ is either the classical Hopf relation (Hopf) 
or the unitary infinitesimal one (u.i.):
\begin{displaymath}
\begin{array}{c||c|c}
 & *&\cdot \\
\hline
\hline
\DD &Hopf & u.i. \\
\hline
\dd & u.i. & Hopf 
\end{array}
\end{displaymath}

Hence, $(2as(V), *, \cdot, \DD, \dd)$ is, at the same time, a free \dda algebra and a cofree \dda coalgebra. 

Observe that $(2as(V), \cdot, *, \dd)$ is also a free \dda bialgebra. In fact the identity on $V$ induces an isomorphism
$$(2as(V), *, \cdot, \DD, \dd)\cong (2as(V), \cdot, *, \dd, \DD).$$

Let $(\HH, *, \cdot, \DD, \dd)$ be a vector space which is a 2-associative algebra, a 2-associative coalgebra verifying the 4 relations of the above tableau. It is not difficult to prove that, if $\HH$ is connected both for $\DD$ and $\dd$, then there is an isomorphism $\HH \cong 2as(\Prim \HH)$ where $\Prim \HH = \Ker \DDo \cap \Ker \overline {\dd}$.

\subsection{The operad $\Bi$}\label{operadBinfinity} 
To unravel the structure of the operad $\Bi$, that is to describe explicitly the free $\Bi$-algebra $\Bi(V)$, we will use the idempotent 
$e:2as(V) \to 2as(V)$ constructed in Proposition \ref{idempotent}. Recall that, a priori, $\Bi(V)= \Bi(n)\t_{S_n} V^{\t n}$ , where the 
right $S_n$-module $ \Bi(n)$  is the space of $n$-ary operations. By Corollary \ref{cor:free2as} $2as(V)$ admits the following decomposition:
$$ 2as(V) = K1 \oplus V \oplus \bigoplus_{n\geq 2} K[\TT_n^*]\t V^{\t n} \oplus \bigoplus_{n\geq 2} K[\TT_n^{\bullet}]\t V^{\t n}.$$
We denote by $t^*$, resp. $t^{\bullet}$, the image of the decorated tree $t\in  K[\TT_n]\t V^{\t n}$ in the relevant component of $2as(V)$. Under this identification, the idempotent $e$ has the following properties:
$$\displaylines{
e(1)= 0,  \\ e(v) = v,\\
e(t^*)= t^* + \oo (t)^{\bullet}\ ,\\
e(t^{\bullet})=0,
}$$
where $t\in  K[\TT_n]\t V^{\t n}, n\geq 2,$ is a decorated tree and where $\oo:  K[\TT_n]\t V^{\t n} \to  K[\TT_n]\t V^{\t n}$ is a 
functorial map in $V$ determined by $e$.

Let us denote by $\cc_{p,q}  \in \TT_{p+q}$ the tree which is the grafting of the $p$-corolla with the $q$-corolla. So $\cc_{pq} $ has 
$p+q$ leaves and 3 vertices:
$$\arbregammapq$$

\begin{thm}\label{thm:freeB} The free  $\Bi$-algebra over the vector space $V$ is 
$$\Bi(V)\cong  \bigoplus_{n\geq 1} K[\TT_n]\t V^{\t n},$$
where $\TT_n$ is the set of planar rooted trees.
Under this isomorphism the operation $M_{pq}\in \Bi(p+q)$ corresponds to $(\cc_{pq}\t 1_{p+q})\in K[\TT_{p+q}]\t K[S_{p+q}]$. The composition of operations in 
$\Bi(V)$ is determined by the following equality which holds in $2as(V)$ (in fact in the $*$-component):
$$M_{pq}(t_1\ldots t_p, t_{p+1}\ldots t_{p+q})^* = \big( e(t_1^*)\cdot \ldots \cdot e(t_p^*)\big) * 
 \big( e(t_{p+1}^*)\cdot \ldots \cdot e(t_{p+q}^*)\big)$$
where the $t_i$'s are $V$-decorated trees.
\end{thm}

First we prove the following Lemma.

\begin{lemma}\label{MpqAnde} Let $x_1 , \ldots , x_{p+q}$ be primitive elements in $2as(V)$. Then we have
$$M_{pq}(x_1 \ldots x_{p+q}) = e \big(( x_1\cdot \ldots \cdot x_p)* (x_{p+1}\cdot \ldots \cdot x_{p+q})\big).$$
\end{lemma}

\P We know that $M_{pq}$ of primitive elements is primitive by Proposition  \ref{prop:primitive}, hence we have 
$e \big( M_{pq}(x_1 \ldots x_{p+q})\big) = M_{pq}(x_1 \ldots x_{p+q})$. By formula (9) of \ref{2asAndB} and item b) in Proposition 
\ref{idempotent} we have $e \big( M_{pq}(x_1 \ldots x_{p+q})\big) = e \big( (x_1\cdot \ldots \cdot x_p)* (x_{p+1}\cdot \ldots 
\cdot x_{p+q})\big)$. Whence the result. \hfill $\square$

 \emph{ Proof of Theorem \ref{thm:freeB}.}
Let us consider the following sequence of maps:
$$K[\TTi(V)] = K[\TTi^{*}(V)] \rightarrowtail 2as(V) \stackrel{e}{\to} \Prim 2as(V) \stackrel{\bar \phi }{\to} \Bi(V).$$
From the properties of $e$ recalled above, namely $e(t^*) = t^* + \oo(t)^{\bullet}$, we know that $e$ restricted to $K[\TTi^{*}(V)] $ is injective. Since $e$ is surjective by Proposition \ref{idempotent} and since  $e(t^{\bullet}) =0$, the restriction of $e$ to $K[\TTi^{*}(V)] $ is surjective too. Finally, $\bar \phi$ is an isomorphism by Proposition \ref{isofree}. As a consequence the above composite is an isomorphism and we have $\Bi(n) \cong K[\TT_n]\t K[S_n]$.

Let us denote by $e'$ the idempotent concerning the unital infinitesimal bialgebra $ T^c(\Bi(V))$. From the functoriality of this construction we deduce a commutative diagram:
\[ \xymatrix{
2as(V) \ar[r]^{\phi}\ar[d]_e & T^c(\Bi(V))\ar[d]^{e'}\\
\Prim (2as(V))\ar[r]^{\bar \phi} & \Bi(V) }\]

We compute:
\begin{eqnarray*}
\bar \phi \circ e (\cc_{p,q} ^*; u_1 \ldots u_pv_1\ldots v_q) & = & \bar \phi \circ e ((u_1\cdot \ldots \cdot u_p)*(v_1
\cdot\ldots\cdot v_q))\\
& = & e'\circ \phi((u_1\cdot \ldots \cdot u_p)*(v_1 \cdot\ldots\cdot v_q))\\
& = & e'(u_1\ldots u_p * v_1\ldots v_q) \\
& = & M_{pq}(u_1\ldots u_p , v_1\ldots v_q) ,
\end{eqnarray*}
since $e'$ is the projection on the component $\Bi(V)$ by Proposition \ref{idempotent} item (d) and by Lemma \ref{MpqAnde}. This 
computation shows that, under the isomorphism $\Bi(n) \cong K[\TT_n]\t K[S_n]$, the operation $M_{pq}$ corresponds to the element 
$\cc_{pq}\t 1_{p+q}$. 

Let $t_i, i=1, \ldots , p+q$, be $V$-decorated trees viewed as elements of $\Bi(V)$. We want to identify the element $M_{pq}(t_1 \ldots  
t_{p+q}) \in \Bi(V)$ as a sum of decorated trees. It is sufficient to describe
$M_{pq}(t_1 \ldots t_{p+q})^*$. The image of $M_{pq}(t_1 \ldots t_{p+q})$ in $2as(V)$ is $e\big( M_{pq}(t_1 \ldots t_{p+q})^*\big)$. 

On the other hand the image of $t_i$ is $e(t_i^*)$ and applying the operation $M_{pq}$ of $2as(V)$ gives $M_{pq}\big(e(t_1^*)\ldots 
e(t_{p+q}^*)\big)$. This element is equal to 
$$e\big( e(t_1^*)\cdot \ldots \cdot e(t_p^*)\big) * 
 \big( e(t_{p+1}^*)\cdot \ldots \cdot e(t_{p+q}^*)\big)$$
 by Lemma \ref{MpqAnde}, whence the equality
$$e\big( M_{pq}(t_1 , \ldots , t_{p+q})^*\big) = e\Big( \big( e(t_1^*)\cdot \ldots \cdot e(t_p^*)\big) * 
 \big( e(t_{p+1}^*)\cdot \ldots \cdot e(t_{p+q}^*)\big)\Big).$$
Since both elements $ M_{pq}(t_1 , \ldots , t_{p+q})^*$ and $ \big(e(t_1^*)\cdot \ldots \cdot e(t_p^*)\big) * 
 \big( e(t_{p+1}^*)\cdot \ldots \cdot e(t_{p+q}^*)\big)$ belong to the $*$-component, we get the expected equality. \hfill $\square$

\subsection{Example} Let $u,v,w,x$ be elements in $V$. Using the formula of Theorem \ref{thm:freeB} we get

$M_{11}\big( (\vert; u) , (\arbreA ; vw)\big) = (\arbreBB ; uvw) - (\arbreAB ; uvw + uwv) , $

$M_{11}\big( (\arbreA ; uv), (\vert; w)\big) =  (\arbreBB ; uvw) - (\arbreBA ; uvw + vuw) , $

\[\begin{array}{l}
M_{11}\big( (\arbreA ; uv), (\arbreA; wx)\big) =(\arbreBBB ; uvwx) - (\arbreABB ;(uv+vu)wx)\qquad \\
\hfill -(\arbreBBA ; uv(wx+xw)) + (\arbreACA ; (uv+vu)(wx+xw)) , 
\end{array}\]

$M_{12}\big( (\arbreA ; uv), (\vert; w) (\vert ;x)\big) = (\arbreBBA ; uvwx) - (\arbreACA ;uvwx+vuwx). $

\subsection{The operations of $\Bi$}
Any planar tree $t\in \TT_n$ determines an $n$-ary operation $M(t)$ in the $\Bi$-operad: $M(t)(v_1\ldots v_n ) = (t; v_1\ldots v_n ) $. 
If $t= \cc_{pq}$, then we know that $M(\cc_{pq})= M_{pq}$ is a generating operation. If not, then $M(t)$ is the composite of the generating 
operations. For instance
\begin{eqnarray*}M(\arbreBB)(uvw) = M_{12}(u, vw+wv) + M_{11}(u, M_{11}(v,w))\\
\hfill = M_{21}(uv+vu, w) + M_{11}(M_{11}(u,v), w).
\end{eqnarray*}

As we see immediately from this example there is no unique way of expressing $M(t)$ in terms of the $M_{pq}$'s because of the relations 
$\RRR_{ijk}$. Here is a recursive algorithm to obtain a formula. Let $t=t^1\vee \ldots \vee t^r\in \TT_n$ be a tree whose root vertex has 
valence $r$. The element 
$$ M_{1\ r-1}\big( (t^1; v_1 \ldots), (t^2; \ldots ), \ldots , (t^r ; \ldots )\big)$$
is of the form $(t^1\vee \ldots \vee t^r; v_1 \ldots) + {\rm other \ terms}$. One can show that all the other terms involve trees whose 
valence is strictly less than $r$. So we can compute $M(t)$ recursively.

\subsection{Example} {the shuffle bialgebra} {\rm The shuffle bialgebra is 
a \dda bialgebra $(T^{sh}(V),  {\scriptstyle \sqcup\!\sqcup}, \cdot, \DD)$ where $ {\scriptstyle \sqcup\!\sqcup}$
is the shuffle product, $\cdot$ the concatenation product and $\DD$ the 
deconcatenation coproduct. The primitive space is $V$.
Its $\Bi$-structure is trivial ($M_{pq}=0$ except for $(p,q)=(1,0)$ and 
$(0,1)$). So there are natural maps
$$2as(V)\twoheadrightarrow U2(V) {\buildrel {\cong} \over \longrightarrow} T^{sh}(V)\ 
.$$
Let us describe explicitly the restriction of the composite $\theta$ to 
the multilinear part of degree $n$ when $V=K$, that is
$$\theta_n :  K[{\TT}^*_n\cup {\TT}^{\bullet}_n]\longrightarrow K.$$
Let $t = t_1\vee \ldots \vee t_k $ be a planar tree whose root valence 
is $k$, and let $p_i$ be the degree of the tree $t_i$.
Since $ {\scriptstyle \sqcup\!\sqcup}$ is the shuffle product and $\cdot $ the concatenation product 
it comes
\begin{eqnarray*}
\theta_n(t^{\cdot}) &=& \theta_n(t^*_1)\cdots \theta_n(t^*_k)\\
\theta_n(t^*)&=& \frac{n!}{p_1!\cdots p_k!}  
\theta_{p_1}(t^{\cdot}_1)\cdots \theta_{p_k}(t^{\cdot}_k)
\end{eqnarray*}
For instance in low dimension $\theta(t)$ is given by:
{\small
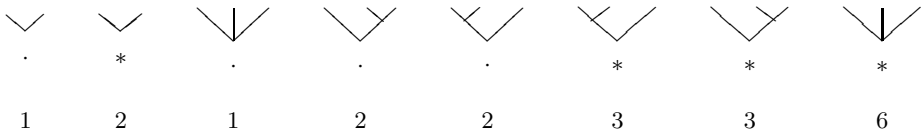
\[ \begin{array}{cccccccccc}
 \arbreAdot & \arbreAe &  \arbreBBdot  &  \arbreABdot &  \arbreBAdot &  \arbreABe &  \arbreBAe &   \arbreBBe \\
&&&&&&&\\
1&2&1&2&2&3&3&6\\
\end{array}\]
}

\section{Dipterous algebras}\label{S:8}

The notion of dipterous algebra is a deformation of the notion of 2-associative algebra. It is, technically, slightly more complicated 
to handle because of the unit. However it has the advantage of containing the case where $\Bi$-algebras are replaced by brace algebras 
(and $2as$-algebras by dendriform algebras). The free dipterous algebra is closely related to the Connes-Kreimer Hopf algebra of rooted 
trees.  Since the proofs of the results in this section (announced in \cite {LR2}) are similar to the 2-associative case, we omit them.

\subsection{Motivation}Let $\HH= (T^c(V), *, \DD)$ be a cofree 
bialgebra. We define inductively a binary operation
$$\d :T^c(V)\t \To^c(V) \to \To^c(V)$$
 by the following formulas
\begin{eqnarray*}
\theta \d v &:=& \theta  v \ ,\\
\theta \d (\oo v) &:=& (\theta * \oo) v\ ,
\end{eqnarray*}
for $\oo\in V^{\t p}, \theta\in V^{\t q}, v\in V$. Observe that $1 v= 
v$, and so $1\d v = v$. The product $\d$ is extended
to $\To^c(V)\t K\,1$ by $\oo\d 1:= 0$, but $1\d 1$ is not defined.

One can check that the following formula holds
$$(x* y)\d z = x\d (y\d z) \eqno (dipt)$$
provided that $y$ and $z$ are not both in $K\, 1=T^c(V)_0$.

\subsection{Dipterous algebra} By definition a \emph {dipterous algebra} is a vector space $A$ equipped with two binary operations 
denoted $*$ and $\d$ satisfying the two relations
\begin{eqnarray*}
(x*y)*z  &:=& x*(y*z) \ ,\\
(x*y)\d z  &:=& x\d (y\d z) \ .
\end{eqnarray*}

So a dipterous algebra is an associative algebra equipped with an extra left action on itself. 

A \emph {unital dipterous algebra} $A$ is a vector space  $A = K 1 \oplus \bar A$ such that $\bar A$ is a dipterous algebra as above, 
where we have extended the two products as follows:
$$ 1*a= a,  \  a*1= a, \quad  \mathrm{ for\ any} \ a\in A .$$
$$  1\d a = a,\ a\d 1 = 0\quad  \mathrm{ for\ any} \ a\in \bar A .$$
Observe that $1\d 1$ is not defined. We denote by ${Dipt}$-alg the category of unital dipterous algebras.

If  $A = K 1 \oplus \bar A$ and  $B = K 1 \oplus \bar B$ are two unital dipterous algebras, then we define a unital dipterous algebra 
structure on their tensor product $A\t B$ as follows:
\begin{eqnarray*}
(a\t b) * (a'\t b') &:=& (a*a') \t (b*b') , \\
(a\t b) \d (a'\t b') &:=& (a*a') \t (b\d b') ,\  \mathrm{if}\ b\t b'\neq 1\t 1 , \\
(a\t 1) \d (a'\t 1) &:=& (a\d a') \t 1 
\end{eqnarray*}
for $a,a'\in A$ and $b,b'\in B$.
\subsection{Dipterous bialgebra}\label{dipterousbialgebra}
By definition a \emph{dipterous bialgebra} $\HH$ is a unital dipterous algebra $(\HH, *, \d)$ equipped with a counital coassociative 
cooperation $\DD$ which satisfies the following compatibility relation: 

$\DD : \HH \to \HH \t \HH$ is a morphism of unital dipterous algebras.
\\

Connectedness and primitive part are defined as in the classical case, cf.~\ref{connectedness} and \ref{Hopf}.

\subsection{Dipterous bialgebra and $\Bi$-algebra}\label{dipterousalgebraAndB} By the same argument as in Proposition 
\ref{infinit-2as} we can show that there is a functor 
$$
(-)_{\Bi} : \{Dipt\mathrm{-alg}\} \to \{\Bi\mathrm{-alg} \}\ .
$$
For instance, defining the operation $\g$ through $x*y = x\g y + x \d y$ , we get:
\begin{eqnarray*}
M_{11}(u,v)  &:=&  u\g v - v\d u\ ,\\
M_{12}(u,vw)  &:=&  u\g (v\d w) - v\d (u\g w) + (v\g w)\d u\ ,\\
M_{21}(uv,w)  &:=& (u\d v)\g w -u\d (v\g w) \ .
\end{eqnarray*}
As in the the $2as$ case one can show that the primitive part of a dipterous bialgebra is a $\Bi$-algebra.

The functor  $(-)_{\Bi}$ has a left adjoint:
$$UD : \{ \Bi\mathrm{-alg}\} \to \{Dipt\mathrm{-alg}\}$$
and $UD(R)$ is a quotient of the free unital dipterous algebra $Dipt(R)$ over the vector space $R$.

For any vector space $V$ there is a unique dipterous homomorphism
$$\DD : Dipt(V) \to Dipt(V)\t Dipt(V)$$
which sends 1 to 1 and $v\in V$ to $v\t 1 + 1 \t v$. It is clearly counital and coassociative, therefore $Dipt(V)$ is a dipterous bialgebra. As a consequence, so is $UD(R)$ for any $\Bi$-algebra $R$.

The free dipterous algebra $Dipt(V)$ admits a description in terms of planar trees similar to the free 2-associative algebra.

\begin{thm}\label{thm:dipterous} If $\HH$ is a (classical) bialgebra over the field $K$, 
then the following are equivalent:

\noindent (a) $\HH$ is a connected dipterous bialgebra,

\noindent (b) $\HH$ is isomorphic to $UD(\Prim \HH)$ as a dipterous bialgebra,

\noindent (c) $\HH$ is cofree among connected coalgebras.
\end{thm}
{}$\hfill \square$

\subsection{Dendriform and brace algebras}\label{dendriform} Let $(A, *, \d)$ be a dipterous algebra, and let $\g$ be the operation defined by the identity $ x*y = x\g y + x\d y$\ . If the operations $\g$ and $\d$ satisfy the relation
$$(x\d y)\g z = x\d (y\g z)$$
then, not only $M_{21}=0$ (cf. \ref{dipterousalgebraAndB}), but all  the operations $M_{pq}$ are 0 for $p\geq 2$. Hence the $\Bi$-algebra associated to the dipterous algebra $A$ is in fact a brace algebra (cf. \ref{examples}).

A dipterous algebra which satisfies the above condition is a \emph {dendriform algebra}, cf. \cite{Lod1}. Equivalently it can be defined as a vector space $A$ equipped with two operations $\g$ and $\d$ satisfying the relations
\begin{displaymath}
\left\{\begin{array}{rcl}
(x\g y)\g z &=& x\g (y* z)   , \\
(x\d y)\g z &=& x\d (y\g z) , \\
(x* y)\d z &=& x\d (y\d z), 
\end{array}
\right.
\end{displaymath}
where $ x*y = x\g y + x\d y $ .

The results of  \cite{R1} and \cite{R2} can be summarized as follows.

\begin{thm}\label{thm:dendriform} If $\HH$ is a dendriform bialgebra over the field $K$, 
then the following are equivalent:

\noindent (a) $\HH$ is connected,

\noindent (b) $\HH$ is isomorphic to $Ud(\Prim \HH)$ as a dendriform bialgebra,

\noindent (c) $\HH$ is cofree among connected coalgebras.
\end{thm}

Here $Ud : \{ Brace{\rm -alg}\} \to \{Dend\mathrm{-alg}\}$ is the left adjoint of the restriction of $(-)_{\Bi}$  to  dendriform algebras.

\subsection{Comparison of Hopf algebras of trees}
As an associative algebra the free dendriform algebra on one generator $Dend(K)$ is a tensor algebra on the planar binary trees (cf. \cite{LR1}). The Connes-Kreimer Hopf algebra $\HH_{CK}$ is the symmetric algebra on (non-planar) rooted trees. Forgetting planarity and symmetrizing gives a surjection of Hopf algebras $Dend(K) \twoheadrightarrow \HH_{CK}$ (cf. for instance \cite{H}).

Since a dendriform algebra is a particular case of dipterous algebra, there is a morphism of dipterous algebras (hence of Hopf algebras):
$$Dipt(K) \to Dend (K)$$
Since $Dipt(K)$ is a tensor algebra over the planar trees, one can describe this map explicitly in terms of trees.
{}$\hfill \square$

\section{Conclusion}\label{S:8}
We compare several variations
of the Cartier-Milnor-Moore theorem.
\\

Summarizing these results, we see that in each 
case three operads are involved: the first one for the
coalgebra, the second one for the algebra, the third one for the 
primitive elements:
\[ \begin{array}{c | c | c | c}
  & \hbox{coalgebra} & \hbox{algebra}& \hbox{primitive}\\
\hline &&&\\
\hbox{Hopf-Borel} & Com & Com & Vect \\
\hbox{Cartier-Milnor-Moore} & Com & As & Lie \\
\hline &&&\\
\hbox{Theorem \ref{thm:infinitesimal}} & As & As & Vect\\
\hbox{Theorem \ref{thm:main}} & As & 2as & \Bi \\
\hline&&& \\
\hbox{M.~Ronco} & As & Dend & brace \\
\hbox{Theorem \ref{thm:dipterous}}& As & Dipt &  \Bi 
\end{array}\]
\\

For each one of these triples of operads there is a Cartier-Milnor-Moore type theorem and a Poincar\'e-Birkhoff-Witt type theorem. 
They suggest the existence of a general result that has to be formulated in terms of operads and cooperads, or,
better, in terms of \emph {props} since the compatibility relation between operations and cooperations is crucial (cf.~\cite{Lod2}). 
The case $(As, 2as, \Bi)$ handled in this paper could serve as a toy-model for this generalization.
\\ 

\noindent {\bf Acknowledgement.} {We thank C.~Brouder, F.~Goichot, E.~Hoffbeck for comments on the first version of this paper. 
Special thanks to Muriel Livernet for her careful reading and clever comments on this paper. This work was partially supported by  the cooperation project France-Argentina ECOS-SeCyt (\# A01E04).}

%\thebibliography

\end{document}